\documentclass[a4paper, 12pt]{article}
\usepackage{makeidx}
\usepackage{enumerate, theorem}
\usepackage{amsmath, amsfonts, amssymb}
\usepackage[height=22.5cm, width=15cm]{geometry}
\usepackage[all]{xy}
\usepackage{mathrsfs, yfonts}
\usepackage{graphicx}
\DeclareMathOperator{\C}{\mathbb{C}}
\newcommand{\parag}[1]{\paragraph{\sc{#1.}} }

\DeclareMathOperator{\Sym}{Sym}

\newtheorem{thm}{Theorem}[subsection]
\newtheorem{defn}[thm]{Definition}
\newtheorem{cor}[thm]{Corollary}
\newtheorem{prop}[thm]{Proposition}
\newtheorem{lemma}[thm]{Lemma}

\begin{document}

\title{On symmetric partial differential operators }

\author{Daniel Barlet\footnote{Institut Elie Cartan, G\'eom\`{e}trie,\newline
Universit\'e de Lorraine, CNRS UMR 7502   and  Institut Universitaire de France.}.}

\maketitle

\parag{Abstract} Let $s_{1}, \dots, s_{k}$ be the elementary  symmetric functions of the complex  variables $x_{1}, \dots, x_{k}$. We say that $F \in \C[s_{1}, \dots, s_{k}]$ is a {trace function} if their exists $f \in \C[z]$ such that $F(s_{1}, \dots, s_{k}] = \sum_{j=1}^{k} f(x_{j})$ for all $s \in \C^{k}$. We give an explicit finite family of second order differential operators in the Weyl algebra $W_{2}:= \C[s_{1}, \dots, s_{k}]\langle \frac{\partial}{\partial s_{1}}, \dots, \frac{\partial}{\partial s_{k}}\rangle $ which generates the left ideal in $W_{2}$ of partial differential operators killing all trace functions. The proof uses a theorem for symmetric differential operators analogous to the usual symmetric functions theorem and the corresponding map for symbols. As a corollary, we obtain for each integer $k$  a  holonomic system which is a quotient of $W_{2}$ by an explicit left ideal whose local solutions are linear combinations of the branches of the multivalued root of the universal equation of degree $k$: $z^{k} + \sum_{h=1}^{k} (-1)^{h}.s_{h}.z^{k-h} = 0$.
\tableofcontents

\section{Introduction}
Let $x_{1}, \dots, x_{k}$ be the roots of the monic polynomial $P_{s}(z) := \sum_{h=0}^{k} (-1)^{h}.\sigma_{h}.z^{k-h}$. For any entire function $f : \C \to \C$ define the trace $T(f)$ of $f$ as the entire function on $\C^{k}$ define by $T(f)(\sigma) := \sum_{j=1}^{k} f(x_{j})$. The aim of this paper is to construct an explicit finite set of second order differential operators in the Weyl algebra
 $$W_{2}:= \C[\sigma_{1}, \dots, \sigma_{k}]\langle \partial_{\sigma_{1}}, \dots, \partial_{\sigma_{k}}\rangle$$ 
 which annihilate any trace functions and such that they give a characterization of  entire functions of this type.\\
 In fact we shall prove more: our explicit finite set of second order differential operators will generate the left ideal in $W_{2}$ of all differential operators killing all trace functions.\\
As it is clear that a trace function view as a function on $x_{1}, \dots, x_{k}$  is killed by any elementary symmetric functions of order at least two of $\partial_{x_{1}}, \dots, \partial_{x_{k}}$ our first approach is to write these symmetric differential operators as elements in $W_{2}$. But one can see directly from the cases $k = 2, 3$ that this is not obvious and even the question of the existence of such operators is not clear. So we are lead to prove an analogous result to the usual symmetric functions theorem for symmetric differential operators (see the theorem  \ref{invariant diff.}).\\
After this result which gives the existence, it appears that explicit computation of the elements $\Sigma_{h} \in W_{2}$ corresponding to the elementary symmetric functions of order   at least two of $\partial_{x_{1}}, \dots, \partial_{x_{k}}$ is out of the range of what  can be computed  explicitly (see the formula of $\Sigma_{3}$ for $k = 3$). Then, using an integral formula to compute trace functions in the spirit of Lisbon integrals, see  [B-MF 19], I  find an explicit finite family of second order differential operators in $W_{2}$ which annihilates any trace function. Then the aim of the second part of the paper is devoted to the proof that any entire function of $\sigma_{1}, \dots, \sigma_{k}$ which is killed by this family is a trace function.\\
To prove this result we have to understand how the symbols behave via the theorem of  symmetric differential operators which is proved in section 3, in order to  describe the characteristic variety of the system given by our explicit family of second order differential operators and to compare with the characteristic variety of the left ideal in $W_{2}$ of all differential operators which are killing trace functions.\\
This allows us to show that our second order differential operators generate the left ideal in $W_{2}$ of {\bf all} differential operators which kill  all the trace functions.\\
An easy corollary of this result is to obtain explicit generators of a left ideal in $W_{2}$ of the differential operators which annihilate the branches of the multivalued root of the universal equation of degree $k$ : $z^{k} + \sum_{h=1}^{k} (-1)^{h}.\sigma_{h}.z^{k-h} = 0$. This system is  holonomic and its $0-$th sheaf of solutions is the $\C-$constructible sheaf generates by local branches of  the multivalued root of the universal equation of degree $k$.\\
An other interesting consequence of our result is that we can also find an explicit finite family of second order differential operators in $W_{2}$ which generates the ideal of differential operators in $W_{2}$ which kills all the ``trace forms'' defined by
$$ T\big(f(z).dz \wedge d\sigma_{1}\wedge \dots \wedge d\sigma_{k-1}\big)(\sigma) := \sum_{j=1}^{k} \frac{f(x_{j})}{P'_{\sigma}(x_{j})}.d\sigma_{1}\wedge \dots \wedge d\sigma_{k} $$
where $f \in \mathcal{O}(\C)$ and where $\sigma_{1}, \dots, \sigma_{k}$ are the symetric functions of $x_{1}, \dots, x_{k}$. \\
We conclude this article by a description of symmetric derivations.

\section{A simple problem}

\subsection{}

 To  each entire function $f : \C \to \C$ we associate the entire function $T(f)$  on $\C^{k}$ defined by $T(f)(x_{1}, \dots, x_{k}) := \sum_{j=1}^{k} \ f(x_{j}) $. We shall call it the {\bf trace} of $f$. Of course, this is a symmetric function in $x_{1}, \dots, x_{k}$ and so the trace of $f$ may be seen as an entire function in the variables $\sigma_{1}, \dots, \sigma_{k}$ defined by the elementary symmetric functions of $x_{1}, \dots, x_{k}$. 

\parag{Example} Any Newton symmetric function $N_{m}(x_{1}, \dots, x_{k}), m \in \mathbb{N}$, considered as a polynomial in $\sigma_{1}, \dots, \sigma_{k}$ is a trace function.\\

If the Taylor expansion of $f$ at $z = 0$ is given by $f(z) = \sum_{m=0}^{\infty} a_{m}.z^{m}$ then we have $T(f)(\sigma) = \sum_{m=0}^{\infty} a_{m}.N_{m}(\sigma)$ where $N_{m}(\sigma)$ is the $m-$th Newton symmetric polynomial corresponding to the roots of the monic polynomial
 $$P_{\sigma}(z) := \prod_{j=1}^{k} (z - x_{j}) = \sum_{h=0}^{k} (-1)^{h}.\sigma_{h}.z^{k-h}$$
 with the convention $\sigma_{0} \equiv1$. \\
\parag{Remark} The  uniform convergence on any compact set of $\C^{k}$ for the expansion of $T(f)(\sigma)$ in quasi-homogeneous polynomials given above, allows to apply any partial derivative term-wise  to such an expansion. So a trace function is an uniform limit on compact sets of linear combinations (with constant coefficients) of the Newton polynomials.\\

The following lemma gives a characterization for an  entire symmetric function on $\C^{k}$ to be a trace function as solutions of a very simple order two differential system.

\begin{lemma}\label{easy 1}
For an entire function $F(x_{1}, \dots, x_{k})$ which is symmetric in $x_{1}, \dots, x_{k}$, a necessary and sufficient condition for the existence of an entire function $f$ on $\C$ such that $F = T(f)$ is given by the following differential system
\begin{equation}
 \frac{\partial^{2}F}{\partial x_{i}\partial x_{j}} = 0 \quad \forall (i, j) \in [1, k]^{2} \quad {\rm with} \quad  i\not= j.
 \end{equation}
\end{lemma}

\parag{Proof} The necessary assumption is obvious. So consider a symmetric entire function solution of the system $(1)$. The partial derivative $\frac{\partial F}{\partial x_{1}}$ does not depends on $x_{2}, \dots, x_{k}$, so we have $\frac{\partial F}{\partial x_{1}}(x) = g(x_{1})$ for some entire function $g$ on $\C$. By symmetry, we conclude that
 $$\frac{\partial F}{\partial x_{j}}(x) = g(x_{j})$$
 for each $j \in [1, k]$. If $f$ is a primitive of $g$ we conclude that $F - T(f)$ is a constant function on $\C^{k}$, so adding a suitable constant to our initial choice for $f$ allows to conclude.$\hfill \blacksquare$
 
 \parag{Remark} It is an easy exercice to prove that the ideal generated by the symbols of the oprerators $\frac{\partial^{2}}{\partial x_{i}\partial x_{j}}$ for $i \not= j$ is reduced. This implies that the left ideal of  
  $\C[x_{1}, \dots, x_{k}]\langle \frac{\partial}{\partial x_{1}}, \dots, \frac{\partial}{\partial x_{k}}\rangle$  which is the annihilator of any trace function is generated by the operators $\frac{\partial^{2}}{\partial x_{i}\partial x_{j}}$ for $i \not= j$.
 
 \subsection{}

Now we want to obtain an analogous characterization of trace functions but with a system given by symmetric differential operators.\\

So introduce in the Weyl algebra $W_{1} := \C[x_{1}, \dots, x_{k}]\langle \frac{\partial}{\partial x_{1}}, \dots, \frac{\partial}{\partial x_{k}}\rangle$ the elementary symmetric ``functions'' $S_{1}, \dots, S_{k}$  of the derivations $\frac{\partial}{\partial x_{j}}, j \in [1, k]$. They are in the sub-algebra $W_{1}^{\mathfrak{S}_{k}}$ of symmetric differential operators. Then we have the following result.

\begin{prop}\label{easy 2}
Let $F$ be an entire function on $\C^{k}$ which is symmetric in $x_{1}, \dots, x_{k}$. Then there exists an entire function $f$ on $\C$ such that $F = T(f)$ if and only if $F$ satisfies the system of partial differential equations
\begin{equation}
S_{h}(F) = 0 \quad \forall h \in [2, k] 
\end{equation}
\end{prop}

\bigskip

The proof will use several lemmas.

\begin{defn}\label{filt}
Let $F_{h}$ be the sub-vector space of \  $T_{k} := \C[x_{1}, \dots, x_{k}]$ generated by the monomials of the type 
$ x_{i_{1}}^{\alpha_{1}}...x_{i_{h}}^{\alpha_{h}}$ where $1 \leq i_{1} < \dots < i_{h}\leq k$ and $\alpha \in \mathbb{N}^{h}$.
\end{defn}

Note that $F_{h}$ is stable by the action of $\mathfrak{S}_{k}$ on $T_{k}$.

\begin{lemma}\label{$F_{k-1}$}
The kernel of the differential operator  $S_{k}$ acting on $T_{k}$  is $F_{k-1}$.
\end{lemma}

\parag{Proof} We shall make an induction on $k \geq 1$. For $k = 1$ we have $F_{k-1} = F_{0} =  \C$ (by convention $x^{\alpha} \equiv 1$ for $\alpha = 0$) and the assertion is clear.\\
Let asume that $k \geq 2$ and the lemma proved for $k-1$. As it is clear that $F_{k-1}$ is contained in $Ker(S_{k})$, consider $P \in Ker(S_{k})$ and write
$$ P = \sum_{p=0}^{d} \ a_{p}.x_{k}^{p} $$
where $a_{p}$ is in $T_{k-1}$ for each $p$. Now $\frac{\partial^{k-1} P}{\partial x_{1} \dots \partial x_{k-1}}$ is independant of $x_{k}$ so belongs to $T_{k-1}$. Then we obtain that $S_{k-1}(a_{p}) = 0$ for each $p \geq 1$ and the induction hypothesis implies that $a_{p}$ is in $F_{k-2} \subset T_{k-1}$ for each $p \geq 1$. So this implies that $a_{p}.x_{k}^{p}$ belongs to $F_{k-1} \subset T_{k}$ for each $p \geq 1$.  As, by definition, $a_{0}$ is in $T_{k-1} \subset F_{k-1}$ we conclude that  $P$ is in $F_{k-1}$.$\hfill \blacksquare$\\

Fix two  integers $1 \leq h \leq k$ and consider the symmetrization map
  $$\sigma_{h, k}: \C[x_{1}, \dots, x_{h}] \to \C[x_{1}, \dots, x_{k}]^{\mathfrak{S}_{k}}$$
 given by
$$ \sigma_{h, k}(P)(x_{1}, \dots, x_{k}) := \frac{(k-h)!}{k!}.\sum_{i \in \mathfrak{I}(h, k)} \ P(x_{i_{1}}, \dots, x_{i_{h}}) $$
where $\mathfrak{I}(h, k)$ is the set of injections of $\{1, 2, \dots, h\}$ into  $\{1, 2, \dots, k\}$.

\begin{lemma}\label{little}
For any integer $h \leq k$ we have $F_{h} \,\cap\,  \C[x_{1}, \dots, x_{k}]^{\mathfrak{S}_{k}} = \sigma_{h, k}\big( \C[x_{1}, \dots, x_{h}]^{\mathfrak{S}_{h}}\big)$.
\end{lemma}

\parag{Proof} This is obvious because $\sigma_{h, k}(P) = \sigma_{h, k}( \sigma_{h, h}(P))$.$\hfill \blacksquare$\\

\begin{lemma}\label{small}
The map \ $\sigma_{h, k}$ is injective on $ \C[x_{1}, \dots, x_{h}] ^{\mathfrak{S}_{h}}$.
\end{lemma}

\parag{Proof} We shall make an induction on $h \in [0, k]$. Assume that $\sigma_{h-1, k}$ is injective (this is obvious for $h= 1$) and consider $P \in \C[x_{1}, \dots, x_{h}]^{\mathfrak{S}_{h}}$ such that $\sigma_{h, k}(P) = 0$. Then $\sigma_{h, k}(P)(x_{1}, \dots, x_{h}, x_{h}, \dots, x_{h}) = ˆ$ and this implies that
$$ (k-h+1).P(x_{1}, \dots, x_{h}) + Q(x_{1}, \dots, x_{h}) = 0 $$
where $Q$ is in $F_{h-1} \subset T_{h}$. Apply $\sigma_{h, k}$ to this equality gives that $P$ belongs to $F_{h-1}\cap \C[x_{1}, \dots, x_{h}]^{\mathfrak{S}_{h}}$, and, thanks to the lemma \ref{little} we obtain that $P = \sigma_{h-1, h}(\tilde{Q})$ where $\tilde{Q}$ is in $\C[x_{1}, \dots, x_{h-1}]^{\mathfrak{S}_{h-1}}$. So we have
$$ 0 = \sigma_{h, k}(P) = \sigma_{h, k}(\sigma_{h-1, h}(\tilde{Q}) = \sigma_{h-1, k}(\tilde{Q}). $$
The induction hypothesis gives $\tilde{Q} = 0$ and so $P = 0$.$\hfill  \blacksquare$\\

\parag{Proof of the proposition \ref{easy 2}}  Consider first a  $Q \in  \C[x_{1}, \dots, x_{h}]^{\mathfrak{S}_{h}}$ such that $S_{h}(\sigma_{h,k}(Q)) = 0$. As we have
$$h!.S_{h}= \sum_{i \in  \mathfrak{I}(h, k)} \ \frac{\partial^{h}}{\partial x_{i(1)} \dots \partial x_{i(h)}} $$
and 
$$  \frac{\partial^{h}}{\partial x_{i(1)} \dots \partial x_{i(h)}}\big[Q(x_{j(1)} \dots x_{j(h)})\big] = 0 \quad for \   i \not= j \ in \  \mathfrak{I}(h, k)\big/\mathfrak{S}_{h}$$
and equal to $\big[\frac{\partial^{h}}{\partial x_{1}, \dots, \partial x_{h}}(Q)\big](x_{1}, \dots, x_{h})$ for $i = j$ \ in \ $\mathfrak{I}(h, k)\big/\mathfrak{S}_{h}$, where we consider the natural action of $\mathfrak{S}_{h}$ on $\mathfrak{I}(h, k)$ by composition,  we obtain
$$ h!.S_{h}(\sigma_{h, k}(Q))(x_{1} \dots x_{k}) = \sigma_{h, k}\big[\frac{\partial^{h}}{\partial x_{1} \dots \partial x_{h}}(Q)\big](x_{1} \dots x_{k}) = 0 .$$
Now the lemma \ref{small} gives $\frac{\partial^{h}}{\partial x_{1} \dots \partial x_{h}}(Q) = 0$.\\
But the kernel of $\frac{\partial^{h}}{\partial x_{1} \dots \partial x_{h}}$ restricted to $\C[x_{1}, \dots, x_{h}]^{\mathfrak{S}_{h}}$ is $F_{h-1} \cap \C[x_{1}, \dots, x_{h}]^{\mathfrak{S}_{h}}$ thanks to lemma \ref{$F_{k-1}$}, so $Q$ is in $F_{h-1} \cap \C[x_{1}, \dots, x_{h}]^{\mathfrak{S}_{h}}$. Then
there exists, thanks to the lemma \ref{little}, $R \in \C[x_{1}, \dots, x_{h-1}]^{\mathfrak{S}_{h-1}}$ such that
 $$Q(x_{1}, \dots, x_{h}) = \sigma_{h-1, h}(R)(x_{1}, \dots, x_{h}) =  \frac{1}{h}. \sum_{j=1}^{h} \ R(x_{1}, \dots, \hat{x}_{j}, \dots, x_{h}).$$
 Then we have
$$ \sigma_{h, k}(Q) = \sigma_{h-1, k}(R) \quad {\rm in} \quad  \C[x_{1}, \dots, x_{k}]^{\mathfrak{S}_{k}}.$$
So we have proved that if  $Q \in  \C[x_{1}, \dots, x_{h}]^{\mathfrak{S}_{h}}$ satisfies $S_{h}(\sigma_{h,k}(Q)) = 0$ there exists $R \in \C[x_{1}, \dots, x_{h-1}]^{\mathfrak{S}_{h-1}}$ such that 
$\sigma_{h-1, k}(R) = \sigma_{h, k}(Q)$.\\
Note that this implies that 
$$ S_{h-1}(\sigma_{h, k}(Q)) =  S_{h-1}(\sigma_{h-1, k}(Q)) = \big((h-1)!\big)^{2}.\sigma_{h-1, k}\big[\frac{\partial^{h-1} R}{\partial x_{1} \dots \partial x_{h-1}}\big] $$
so that  $  S_{h-1}(\sigma_{h, k}(Q)) = 0$ implies $\frac{\partial^{h-1} R}{\partial x_{1} \dots \partial x_{h-1}} = 0$ as $\sigma_{h-1, k}$ is injective.\\
Then if $P$ is solution of $S_{h}(P) = 0$ for each $h \in [2, k]$ we construct by a descending induction a sequence $Q_{h} \in \C[x_{1}, \dots, x_{h}]^{\mathfrak{S}_{h}}$ such that $P = \sigma_{h, k}[Q_{h}]$. So we find finally, thanks to the first part of this proof,  a $Q_{1} \in \C[x_{1}]$ with
 $$P(x_{1}, \dots, x_{k}) = \frac{1}{k}.\sum_{j=1}^{k} Q_{1}(x_{j}) = T(\frac{1}{k}.Q_{1})(x_{1}, \dots, x_{k}).$$

Consider now a symmetric  entire function $F$ such that $S_{h}(F) = 0$ for each $h \in [2, k]$. Let $F = \sum_{\nu \geq 0} P_{\nu}$ be its Taylor expansion at the origin in $\C^{k}$. Each homogeneous polynomial $P_{\nu} $ is symmetric and as the differential operators $S_{h}$ are homogeneous we obtain that $S_{h}(P_{\nu}) \equiv 0$ for each $\nu \geq 0$ and each $h \in [2, k]$ by uniqueness of the Taylor expansion. Then we can find polynomials $R_{\nu} \in \C[x_{1}]$ which are homogeneous of degree $\nu$ such that $\sigma_{1, k}[R_{\nu}] = P_{\nu}$. From the convergence of the Taylor series of $F$, it is then easy to see that the series $\sum_{\nu \geq 0} R_{\nu}$ converges uniformly on compact sets in $\C$ and defines an entire function $f$ such that $T(f) = F$.\\
As the converse is obvious, the proposition is proved.$\hfill \blacksquare$\\

\begin{prop}\label{complement}
Any element in $W_{1}^{\mathfrak{S}_{k}}$ which kills any Newton polynomial 
$$N_{m}(x) := \sum_{j=1}^{k} x_{j}^{m}, \quad  m \in \mathbb{N}$$
  is in the left ideal of  $W_{1}^{\mathfrak{S}_{k}}$ generated by $S_{2}, \dots, S_{k}$.
\end{prop}

The proof will be a simple consequence of the following lemma.

\begin{lemma}\label{simplet 0}
The intersection of the ideal $I$ in $\C[\xi_{1}, \dots, \xi_{k}]$,  generated by $\xi_{i}.\xi_{j}$ for $i \not= j$  with $\C[\xi_{1}, \dots, \xi_{k}]^{\mathfrak{S}_{k}}$ in the ideal generated by the elementary symmetric functions $S_{2}, \dots, S_{k}$  of \ $\xi_{1}, \dots, \xi_{k}$.
\end{lemma}

\parag{Proof} Remark first that $S_{2}, \dots, S_{k}$ are in $I$ and symmetric. Let $P \in \C[\xi_{1}, \dots, \xi_{k}]^{\mathfrak{S}_{k}}$ be in $I$. Then as $P$ is in $\C[S_{1}, \dots, S_{k}]$ we may write
$$ P = Q(S_{1}) \quad {\rm modulo} \ (S_{2}, \dots, S_{k})$$
where $Q$ is in $\C[z]$. Fixing $\xi_{2} = \dots = \xi_{k} = 0$ we obtain, as $P$ is in $I$, $Q(\xi_{1}) = 0$ for any value of $\xi_{1}$. So $Q = 0$ and $P$ belongs to  the ideal generated by $S_{2}, \dots, S_{k}$. Moreover this ideal is obviously prime in $\C[S_{1}, \dots, S_{k}] = \C[\xi_{1}, \dots, \xi_{k}]^{\mathfrak{S}_{k}}$.$\hfill \blacksquare$\\

\parag{Remark} Note that the ideal $I$ is reduced in $\C[\xi_{1}, \dots, \xi_{k}]$ (but not prime).$\hfill \square$\\

\parag{Proof of the proposition \ref{complement}} Let $P$ be in $W_{1}^{\mathfrak{S}_{k}}$ and assume that $P$ kills any Newton polynomial, is not in the left  ideal generated by $S_{2}, \dots , S_{k}$ in $W_{1}^{\mathfrak{S}_{k}}$ and has minimal order for these properties.  Then  the remark following the lemma \ref{easy 1} implies that $P$ is in the left ideal of $W_{1}$ generated by the $\frac{\partial^{2}}{\partial x_{i}\partial x_{j}}, i \not= j$. So the symbol of $P$ is in $I \cap  \C[\xi_{1}, \dots, \xi_{k}]^{\mathfrak{S}_{k}}$. Using the lemma above, there exists $Q$ in the left  ideal generated by $S_{2}, \dots , S_{k}$ in $W_{1}^{\mathfrak{S}_{k}}$ such that $Q$ has the same symbol than $P$. So the order of $P- Q$ is strictly  less than the order of $P$ and its kills any Newton polynomial. Then $P-Q$ is in  left  ideal generated by $S_{2}, \dots , S_{k}$ in $W_{1}^{\mathfrak{S}_{k}}$. Contradiction ! So the lemma is proved.$\hfill \blacksquare$\\

Nevertheless we are not happy with this result because, using the symmetric function theorem, we are in fact working with an entire function $F$ on $\C^{k}$ with coordinates $\sigma_{1}, \dots, \sigma_{k}$ corresponding to the elementary symmetric functions of $x_{1}, \dots, x_{k}$ and we would like to have a system of partial differential operators in the Weyl algebra $W_{2} := \C[\sigma_{1}, \dots, \sigma_{k}]\langle \frac{\partial}{\partial \sigma_{1}}, \dots, \frac{\partial}{\partial \sigma_{k}} \rangle$ characterizing the entire functions $F$ which are trace functions. \\

So the problem is know to rewrite the system $(2)$ of the proposition \ref{easy 2} in term of elements in $W_{2}$. This will be the aim of the next section.\\

\section{The symmetric function theorem for linear differential operators}

It will be convenient to look at the local version in the study of trace functions.  This is the aim of our first paragraph.

\parag{Notations} \begin{itemize}
\item We shall denote by $x := (x_{1}, \dots, x_{k})$ a point in $M := \C^{k}$, by $\sigma := (\sigma_{1}, \dots, \sigma_{k})$ a point in $N := \Sym^{k}(\C) \simeq \C^{k}$ and by $s : M \to N$ the quotient map given by the elementary symmetric functions $x \mapsto (\sigma_{1} = s_{1}(x), \dots, \sigma_{k} = s_{k}(x))$.
\item We shall denote by $\Delta : N \to \C$ the discriminant polynomial. So 
\begin{equation*}
\Delta(s(x)) = \prod_{1\leq i < j \leq k}(x_{i} - x_{j})^{2}. \tag{$\Delta$}
\end{equation*}
\end{itemize}

\subsection{Local trace functions}

In this section we shall consider the hypersurface
\begin{equation*}
 H := \{(\sigma, z) \in N \times \C \ / \  P_{\sigma}[z] = 0 \}. \tag{H}
 \end{equation*}
where $P_{\sigma}[z] := z^{k} + \sum_{h=1}^{k} (-1)^{h}.\sigma_{h}.z^{k-h}$.\\
It is smooth an isomorphic with $\C^{k}$  via the map $(\sigma, z) \mapsto (\sigma_{1}, \dots, \sigma_{k-1}, z)$. We shall denote by $\pi : H \to N$ and $p : H \to \C$, the natural projections.

\begin{defn}\label{local}
Let $V$ be an open set in $N$ and let $F$ be a holomorphic function on $V$. We shall say that $F$ is a {\bf  (global) trace function} on $V$ if there exists a holomorphic function $f$ on $p(\pi^{-1}(V))$ such that
\begin{equation*}
 F(\sigma) = Trace(f)(\sigma) := \sum_{P_{\sigma}[x_{j}] = 0} f(x_{j}), \quad \forall \sigma \in V  \tag{T}
 \end{equation*}
counting multiplicities.\\
 If the holomorphic function $F$ on $V$ is a trace function in a neighborhood of any point in $V$ we shall say that $F$ is a { \bf local trace function} on $V$
\end{defn}

\parag{Examples}\begin{enumerate}
\item For each $m \in \mathbb{N}$ the trace function of the polynomial function  $z \mapsto z^{m}$ is the $m-$th Newton polynomial $N_{m}$ in $\C[\sigma_{1}, \dots, \sigma_{k}]$.
\item Let $Q \in \C[z]$ and fix a point $\sigma^{0}\in N$ such that $P_{\sigma^{0}}$ admits a simple root $x^{0}$. Let $V$  be an open neighborhood of $\sigma^{0}$ on which their exists a holomorphic function $\varphi : V \to \C$ such that $\varphi(\sigma^{0}) = x^{0}$ and $P_{\sigma}[\varphi(\sigma)] = 0$ for each $\sigma \in V$. Then $\varphi$ is a  trace function on $V$. But if $g$ is any holomorphic function on the open set  $\varphi(V) \subset \C$, $\sigma \mapsto g(\varphi(\sigma)) $ is also a trace function on $V$.
\end{enumerate}

\begin{lemma}\label{unique} If $F$ is a (global) trace function on the open set $V$;  then the holomorphic  function $f$ on $p(\pi^{-1}(V))$ such that $F = Trace(f)$ on $V$ is unique up to a locally constant function. So when $p(\pi^{-1}(V))$ is connected, $f$ is unique.
\end{lemma}

\parag{Proof} It is enough to consider a holomorphic function $g$ on $p(\pi^{-1}(V))$ such that  $Trace(g) = 0$ on $V$ and to prove that $g$ is locally constant on $p(\pi^{-1}(V))$.\\
 Fix a point $\sigma^{0} \in V $ such that $\Delta(\sigma^{0}) \not= 0$. Then choose open discs $D_{1}, \dots, D_{k}$ in $\C$ such that
\begin{enumerate}
\item $D_{1}, \dots, D_{k}$ are disjoint.
\item Any $D_{j}$ contains exactly one root of $P_{\sigma^{0}}$ which is its center.
\end{enumerate}
Let $V := q(D_{1}\times \dots\times D_{k})$. Then their exist holomorphic maps $\varphi_{j} : V \to D_{j}$ for $j \in [1,k]$  such that we have for any $\sigma \in V$
$$ P_{\sigma}[z] = \prod_{j=1}^{k}(z - \varphi_{j}(\sigma)), \quad \forall z \in \C.$$
Now fix $x_{2}^{0}, \dots x_{k}^{0}$ respectively in $D_{2}, \dots, D_{k}$. Then $Trace(g)$ at the point $s(x_{1}, x_{2}^{0}, \dots x_{k}^{0})$ vanishes for any $x_{1} \in D_{1}$. This implies that $g$ is constant on $D_{1}$. So $g$ is constant on each $D_{j}, j \in [1, k]$. Then $g$ is locally constant on $p(\pi^{-1}(V))$ as  the complement of a hypersurface in a connected open set is  connected.$\hfill \blacksquare$\\

\begin{cor}\label{local-global}
Let $V$ be a connected open set in $N$ such that $p(\pi^{-1}(V))$ is (connected and)  simply connected. Then  any local trace function $F$ on $V$ is a global trace function on $V$.
\end{cor}

\parag{Proof} Let $(V_{a})_{a \ in A}$ be an open covering of $V$ such that for each $a \in A$ their exists a holomorphic function $f_{a}$ on $p(\pi^{-1}(V_{a}))$ such that $F$ is the trace of $f_{a}$ on $V_{a}$. Then by the previous lemma we know that $f_{a}- f_{b}$ is a locally constant function on $p(\pi^{-1}(V_{a})) \cap p(\pi^{-1}(V_{b}))$. This define a $1-$cocyle on the open covering of $p(\pi^{-1}(V))$ with value in the constant sheaf. As we assume that $p(\pi^{-1}(V))$ is connected and simply connected, up to pass to a finer covering, we may assume that their exist locally constant functions $g_{a}$ on $p(\pi^{-1}(V_{a}))$ with $f_{a} - f_{b} = g_{a} - g_{b}$ on the intersection $p(\pi^{-1}(V_{a})) \cap p(\pi^{-1}(V_{b}))$. Let $\tilde{f}$ be the global holomorphic function on 
$p(\pi^{-1}(V))$ given by $f_{a}- g_{a}$ on the  $p(\pi^{-1}(V_{a}))$
 Then the trace of $\tilde{f}$ is a holomorphic function $\tilde{F}$ on $V$ and $\tilde{F} - F$ is locally constant on $V$. So it is constant equal to some $\gamma \in \C$ on $V$. Then the trace of  function $\tilde{f} - \gamma/k$ is equal to $F$.$\hfill \blacksquare$

\parag{Example} When $V$ is the image by the quotient map $s : M \to N$ of $D^{k}$ where $D$ is an open  disc in $\C$, $V$ is connected and $p(\pi^{-1}(V)) = D$ is simply connected.

\begin{lemma}\label{Newton}
Let $V$ be a connected open set in $N$ such that $p(\pi^{-1}(V))$ is connected and simply connected. Then $F \in \mathcal{O}(V)$ is a trace function if and only if it is a uniform limit on compact sets in $V$ of linear combinations of Newton polynomials.
\end{lemma}

\parag{Proof} As $p(\pi^{-1}(V))$ is a Runge domain, any holomorphic function on it is a uniform limit on compact sets of polynomials.  But  the trace of a polynomial is a linear combination (with constant coefficients) of the Newton polynomials, and the result follows.$\hfill \blacksquare$\\

\subsection{The theorem}

\parag{Notations} \begin{itemize}
\item Let  $e_{1}, \dots, e_{k}$ be then the standard basis of $T_{M} \simeq M\times \C^{k}$ and $\varepsilon_{1}, \dots, \varepsilon_{k}$  the  standard basis of  $T_{N} \simeq N\times \C^{k}$.
\item Let $T_{s}: T_{M} \to s^{-1}(T_{N})$ the tangent map to $s$. We shall identify  a vector in $T_{M, x}$ with a linear form of $T_{M, x}^{*}$ 
 In the same manner we shall identify a vector in $T_{N, \sigma}$ (or in $s^{-1}(T_{N})_{x}$) with a linear form on $T_{N, \sigma}^{*}$ (resp. in $s^{-1}(T_{N}^{*})_{x} \simeq s^{-1}(T_{N})_{x}^{*}$).
With this convention, the symbol of $\frac{\partial}{\partial x_{i}}$ is the vector field $e_{i}$ identified with the function $e_{i}$ on $T_{M}^{*}$ (linear on the fibers), and the symbol of $\frac{\partial}{\partial \sigma_{h}}$ is the vector field $\varepsilon_{h}$ identified with the function $\varepsilon_{h}$ on $T¬^{*}_{N}$ (linear on the fibers).
\item The tangent map $T_{s}$ sends the vector $(x, e_{i})$ to the vector $T_{s}(x, e_{i}) \in T_{N, s(x)}$
$$T_{s}(x, e_{i}) =  (s(x), \sum_{h=1}^{k} \frac{\partial s_{h}}{\partial x_{i}}(x).\varepsilon_{h}).$$
As we have (see the lemma \ref{simplet 2} below)
$$ \frac{\partial s_{h}}{\partial x_{i}}(x) = \sum_{p=0}^{h-1} (-x_{i})^{p}.s_{h-p-1}(x) $$
this gives, for  $(x, e_{i}) \in T_{M, x} , i \in [1, k]$
\begin{equation}
T_{s}(x, e_{i})= (s(x),  \sum_{h=1}^{k}\big(\sum_{p=0}^{h-1} (-x_{i})^{p}.s_{h-p-1}(x)\big).\varepsilon_{h}) .
\end{equation}
\item We shall write $(x, \xi), \xi := \sum_{i=1}^{k} \xi_{i}.e_{i}$, a point in $T_{M, x}$, $(\sigma, \eta), \eta := \sum_{h=1}^{k} \eta_{h}.\varepsilon_{h}^{*}$, a point\footnote{$\varepsilon_{1}^{*}, \dots, \varepsilon_{k}^{*}$ denotes the dual basis of $\varepsilon_{1}, \dots, \varepsilon_{k}$.}  in $T_{N, \sigma}^{*}$ and also $(x, s(x), \eta)$ a point in $s^{-1}(T_{N}^{*})_{x}$.
\item Let now introduce for $h \in [1, k] $ and $\sigma \in N$ the polynomial 
$$\Theta_{h}(z, \sigma) := \sum_{p=0}^{h-1} (-z)^{h-p-1}.\sigma_{p}.$$
Then the formula $(3)$ may be written in the following way
$$ T_{s}(x, \xi) = \sum_{h=1}^{k}\Big(\sum_{i=1}^{k} \Theta_{h}(x_{i}, s(x)).\xi_{i}\Big).\varepsilon_{h} $$
and then the cotangent map $T_{s}^{*} : s^{-1}(T_{N})^{*} \to T_{M}^{*}$ which is defined by the equality
 $$\langle T_{s}^{*}(x, s(x), \eta), \xi \rangle = \langle T_{s}(x, \xi), \eta \rangle $$
 may be written
\begin{equation}
T_{s}^{*}(x, s(x), \eta) = \sum_{i=1}^{k} \Big(\sum_{h=1}^{k} \Theta_{h}(x_{i}, s(x)).\eta_{h}\Big).e_{i}^{*}
\end{equation}
so that the components in the basis $e_{1}^{*}, \dots, e_{k}^{*}$ \ are the numbers $\sum_{h=1}^{k} \Theta_{h}(x_{i}, s(x)).\eta_{h}$.
\item Let $D_{M}$ (resp. $D_{N}$) be the sheaf of holomorphic differential operators\footnote{The reader may consult [Bj] or [Bor] for basic results on $D-$modules.} on $M$ (resp. on $N$).
\item Let $G \subset M \times N$ be the graph of the quotient map $s$, and let $p_{1}: G \to M$ and $p_{2}: G \to N$ the natural projections. Note that $p_{1}$ is an isomorphism and that $p_{2}$ is a proper finite surjective map.
\item The natural action of $\mathfrak{S}_{k}$ on $M = \C^{k}$ defines an action on $M\times N$ (the action on $N$ is, of course, trivial) and also on $G$. For this action $p_{1}$ and $p_{2}$ are equivariant. Then we have an action of $\mathfrak{S}_{k}$ on the sheaves $(p_{2})_{*}(p_{1}^{*}(D_{M}))$ and $(p_{2})_{*}(p_{2}^{*}(D_{N}))$.\\
 We shall denote respectively by $(p_{2})_{*}(p_{1}^{*}(D_{M}))^{\mathfrak{S}_{k}}$ and $(p_{2})_{*}(p_{2}^{*}(D_{N}))^{\mathfrak{S}_{k}}$ the sub-sheaves of $\mathfrak{S}_{k}-$invariant sections of these sheaves. They are $\mathcal{O}_{N}-$sub-modules of  these $\mathcal{O}_{N}-$algebras.
\item We have also an action of $\mathfrak{S}_{k}$ on the sheaf of $\mathcal{O}_{N}-$algebras $s_{*}(D_{M})$ and we shall note $s_{*}(D_{M})^{\mathfrak{S}_{k}}$ the corresponding  $\mathcal{O}_{N}-$sub-modules of the sheaves of $\mathfrak{S}_{k}-$invariant sections.
\item As $s$ is the (holomorphic) quotient map (see [B.M 1] chapter I  theorem 3.1.13) we have a natural isomorphism $\theta : s_{*}(\mathcal{O}_{M})^{\mathfrak{S}_{k}} \to \mathcal{O}_{N}$ whose inverse is simply the pull-back of holomorphic functions.
\end{itemize}

\begin{thm}\label{invariant diff.}
There exists a natural morphism of  sheaves of $\mathcal{O}_{N}-$algebras
$$ \Xi : s_{*}(D_{M})^{\mathfrak{S}_{k}} \to D_{N} $$
such that, for any $P \in  s_{*}(D_{M})^{\mathfrak{S}_{k}}$ and   for any $f \in s_{*}(\mathcal{O}_{M})^{\mathfrak{S}_{k}}$ we have
\begin{equation}
\theta(P[f]) = \Xi(P)[\theta[f]] 
\end{equation}
Moreover, this morphism is injective and it restriction to the Zariski open set
 $$\{\sigma \in N \ / \  \Delta(\sigma) \not= 0 \}$$
  in $N$ is an isomorphism of sheaves of $\mathcal{O}_{N}-$algebras.\\
\end{thm}

The proof of this theorem will use the following simple lemma.

\begin{lemma}\label{simplet 1}
We have  natural isomorphisms of sheaves of $\mathcal{O}_{N}-$algebras
\begin{align*}
& ( p_{2})_{*}(p_{1}^{*}(D_{M})) \to s_{*}(D_{M}) \\
&  ( p_{2})_{*}(p_{1}^{*}(D_{M}))^{\mathfrak{S}_{k}} \to s_{*}(D_{M})^{\mathfrak{S}_{k}}\\
&  (p_{2})_{*}(p_{2}^{*}(D_{N}))^{\mathfrak{S}_{k}} \to D_{N}
\end{align*}
\end{lemma}

\parag{Proof} The first isomorphism is clear as $p_{1}$ is an isomorphism and as we have  the equality $p_{1}^{-1}\circ p_{2} = s$. This implies the second isomorphism by $\mathfrak{S}_{k}-$equivariance.\\
To prove the last isomorphism, remark that we have
 $$(p_{2})_{*}(p_{2}^{*}(D_{N})) \simeq (p_{2})_{*}(\mathcal{O}_{G}) \otimes_{\mathcal{O}_{N}} D_{N}$$
 and so :
 $$(p_{2})_{*}(p_{2}^{*}(D_{N}))^{\mathfrak{S}_{k}}  \simeq (p_{2})_{*}(\mathcal{O}_{G})^{\mathfrak{S}_{k}} \otimes_{\mathcal{O}_{N}} D_{N} ;$$
 as $(p_{2})_{*}(\mathcal{O}_{G})^{\mathfrak{S}_{k}} \simeq \mathcal{O}_{N}$ using $\theta$ and the $\mathfrak{S}_{k}-$equivariant  isomorphism $p_{1}$ which gives $s_{*}(\mathcal{O}_{M}) \simeq (p_{2})_{*}(\mathcal{O}_{G})$.$\hfill \blacksquare$\\
 
 \parag{Proof of the theorem} The tangent map to the map $p_{2}$ may be view as a map of vector bundles  $\tilde{T} : (p_{1})^{*}(T_{M}) \to p_{2}^{*}(T_{N})$ on $G$. It induces a morphism of sheaves of $\mathcal{O}_{G}-$algebras
 $$ \hat{T} : p_{1}^{*}(D_{M}) \to p_{2}^{*}(D_{N}) $$
 as vector fields generate the algebra of differential operators. This morphism is equi-variant for the action of $\mathfrak{S}_{k}$ on $G$ and, thanks to the equivariance of $p_{2}$ we obtain a morphism of $\mathcal{O}_{N}-$algebras
 $$\Xi_{0} : (p_{2})_{*}(p_{1}^{*}(D_{M}))^{\mathfrak{S}_{k}} \to (p_{2})_{*}(p_{2}^{*}(D_{N}))^{\mathfrak{S}_{k}} .$$
 This  gives the definition of  the morphism $\Xi$ \  using the last two isomorphisms of the previous lemma.\\
 To prove the formula $(5)$ it is enough to prove it on the dense Zariski  open set $\{\sigma \in N \ / \  \Delta(\sigma) \not= 0 \}$ in $N$. So consider $D_{1}, \dots, D_{k}$ disjoints open discs in $\C$  and let $U_{0} := D_{1}\times \dots \times D_{k}$ and $\mathcal{U} := s(U_{0})$. Then $s$ induces an isomorphism of $U_{0}$ onto $\mathcal{U}$. If $U_{1} := s^{-1}(\mathcal{U})$, a $\mathfrak{S}_{k}-$invariant function (resp. $\mathfrak{S}_{k}-$invariant differential operator) on $U_{1}$ is simply given by a function (resp. a differential operator) on $U_{0}$ or equivalently by a function (resp. by a differential operator)  in $\mathcal{U}$ via the isomorphism $s_{\vert U_{0}}$. Then $\Xi_{\vert \mathcal{U}}$ is an isomorphism of $s_{*}(D_{M})_{\vert \mathcal{U}}^{\mathfrak{S}_{k}} \simeq s_{*}((D_{M})_{\vert U_{0}}) \simeq (D_{N})_{\vert \mathcal{U}}$ and the formula $(5)$ is satisfied on $\mathcal{U}$ as $\theta$ reduces to the composition of the restriction to $U_{0}$ with  $s$. $\hfill \blacksquare$\\
 
 \begin{cor}\label{Symbol}
In the situation of the previous theorem the map $\Xi$ induces a graded sheaf map between commutative graded $\mathcal{O}_{N}-$algebras :
$$ S\Xi : s_{*}(Gr^{\bullet}D_{M})^{\mathfrak{S}_{k}} \to Gr^{\bullet}D_{N} $$
which sends the symbol of a non zero  symmetric differential operator $P \in s_{*}(D_{M})^{\mathfrak{S}_{k}}$ to the symbol of \  $\Xi(P)$ in $GR^{\bullet}D_{N}$. 
\end{cor}

\parag{Proof} Remark that the morphism $\Xi$ preserves the order of a non zero differential operator. So the morphism $\Xi$ respects the filtrations of the sheaves $s_{*}(D_{M})^{\mathfrak{S}_{k}} $ and $D_{N}$ by the order of the differential operators. Then the morphism $S\Xi$ is obtained from the morphism $\Xi$ by passing to the graded algebras associated to these filtrations.$\hfill \blacksquare$\\

Let us finish this paragraph by the simple lemma used above.

\begin{lemma}\label{simplet 2}
Let $x_{1}, \dots, x_{k}$ be complex numbers and let $s_{1}, \dots, s_{k}$ their elementary symmetric functions. For $j \in [1, k]$ let $s_{h}(j)$ the $h-$th symmetric function of $x_{1}, \dots, \hat{x}_{j}, \dots, x_{k}$ with the convention $s_{0} = s_{0}(j) = 1$. Then we have
 \begin{equation}
\frac{\partial s_{h}}{\partial x_{j}}  = s_{h-1}(j) = \sum_{q=0}^{h-1} \ s_{h-q-1}.(-x_{j})^{q} \qquad \forall \ h \in [1, k] .
 \end{equation}
 Let $\mathfrak{S}_{k}(j)$ be the stabilizer of $j$ in $\mathfrak{S}_{k}$, the permutation group of $\{1, 2, \dots, k\}$. The sub-algebra of $\mathfrak{S}_{k}(j)-$invariant elements in $\C[x_{1}, \dots, x_{k}] $ is equal to the sub-algebra generated by $\C[s_{1}, \dots, s_{k}]$ and  $x_{j}$. Moreover this sub-algebra  is a  free \ $\C[s_{1}, \dots, s_{k}]-$module with basis $1, x_{j}, \dots, x_{j}^{k-1}$.
 \end{lemma}
 
 \parag{Proof} The first equality in $(6)$ is clear. The second equality in $(6)$  is obvious for $h = 1$. Then assume that it is proved for $h-1 \geq 1$. Using the easy equality $s_{h}(j) = s_{h} -x_{j}.s_{h-1}(j)$ we obtain
 $$ s_{h}(j) = s_{h} - x_{j}.\big(\sum_{q=0}^{h-1} \ s_{h-q-1}.(-x_{j})^{q}\big) = \sum_{q=0}^{h} \ s_{h-q}.(-x_{j})^{q} .$$
concluding the induction. The last assertions follows immediately. $\hfill \blacksquare$\\
 
\parag{Remark} The morphism $S\Xi$ of $\mathcal{O}_{N}-$algebras is completely determined by the formula $(3)$.

 \subsection{Examples}

\parag{The case $k = 2$}

For $k = 2$ we obtain that to $\frac{\partial^{2}}{\partial x_{1}\partial x_{2}} $ corresponds via $\Xi$ the differential operator
$$ \Sigma_{2} = \frac{\partial^{2}}{\partial \sigma_{1}^{2}} + \sigma_{1}. \frac{\partial^{2}}{\partial \sigma_{1}\partial \sigma_{2}} + \sigma_{2}.\frac{\partial^{2}}{\partial \sigma_{2}^{2}} + \frac{\partial}{\partial \sigma_{2}}$$
in $W_{2}$.

\parag{Exercice} Check that $\Sigma_{2}$ kills the Newton polynomials $N_{m}(\sigma), \forall m \in \mathbb{N}$ \\
{\it hint} : use J. Varouchas formula's :
$$N_{m}(x_{1}, \dots, x_{k}) = \sum_{\alpha \in \mathbb{N}^{k}, \vert\vert \alpha \vert\vert = m} (-1)^{m+ \vert\alpha\vert}.m.\frac{(\vert\alpha\vert -1)!}{\alpha!}.\sigma^{\alpha} $$
where $\vert \alpha\vert := \sum_{j=1}^{k} \ \alpha_{j}$ and $\vert\vert\alpha\vert\vert := \sum_{j=1}^{k} \ j.\alpha_{j}$.

\parag{The case $k=3$} Let $\Sigma_{h} = \Xi(S_{h})$, where $S_{h}$ is the $h-$th elementary symmetric function of $\frac{\partial}{\partial x_{1}}, \dots, \frac{\partial}{\partial x_{k}}$.\\
Some tiring computations give
  \begin{align*}
  &  \Sigma_{2} = 3\partial_{1}^{2} + (\sigma_{1}^{2} + \sigma_{2}).\partial_{2}^{2} + \sigma_{1}.\sigma_{3}\partial_{3}^{2} + 4\sigma_{1}.\partial_{1}\partial_{2} + 2\sigma_{2}.\partial_{1}\partial_{3} + (\sigma_{1}.\sigma_{2} + 3\sigma_{3}).\partial_{2}\partial_{3} + 3\partial_{2} + \sigma_{1}.\partial_{3}\\
  & \Sigma_{3} = \partial_{1}^{3} + (\sigma_{1}.\sigma_{2} - \sigma_{3}).\partial_{2}^{3} + \sigma_{3}^{2}.\partial_{3}^{3} + 2\sigma_{1}.\partial_{1}^{2}\partial_{2} + \sigma_{2}.\partial_{1}^{2}\partial_{3} + \\
  & \qquad  + (\sigma_{2} + \sigma_{1}^{2}).\partial_{1}\partial_{2}^{2} + (\sigma_{1}.\sigma_{2} + 3\sigma_{3}).\partial_{1}\partial_{2}\partial_{3} + \sigma_{1}.\sigma_{3}.\partial_{1}\partial_{3}^{2}  + \\
  & \qquad  + (\sigma_{1}.\sigma_{3} + \sigma_{2}^{2}).\partial_{2}^{2}\partial_{3} + 2\sigma_{2}.\sigma_{3}.\partial_{2}\partial_{3}^{2} + 3\partial_{1}\partial_{2} + \sigma_{1}.\partial_{1}.\partial_{3} + \\
  & \qquad  + 3\sigma_{2}.\partial_{2}\partial_{3} + 2\sigma_{1}.\partial_{2}^{2} + 3\sigma_{3}.\partial_{3}^{2} + \partial_{3}
  \end{align*}
  
  Let us check $\Sigma_{3}N_{6}$. As we have
   $$N_{6} =   \sigma_{1}^{6} - 6\sigma_{1}^{4}.\sigma_{2} + 6\sigma_{1}^{3}.\sigma_{3} + 9\sigma_{1}^{2}.\sigma_{2}^{2} -12\sigma_{1}.\sigma_{2}.\sigma_{3} - 2\sigma_{2}^{3} + 3\sigma_{3}^{2}$$
   \begin{align*}
  & \Sigma_{3}N_{6} = (120\sigma_{1}^{3} + (\sigma_{1}.\sigma_{2} -\sigma_{3})(-12) + 2\sigma_{1}.(-72\sigma_{1}^{2} +36\sigma_{2}) + \sigma_{2}.(36\sigma_{1}) + (\sigma_{2} + \sigma_{1}^{2}).(36\sigma_{1}) + \\
  & \qquad + (\sigma_{1}.\sigma_{2} + 3\sigma_{3}).(-12) +  3(-24\sigma_{1}^{3} + 36\sigma_{1}.\sigma_{2} -12\sigma_{3}) + \sigma_{1}.(18\sigma_{1}^{2} - 12\sigma_{2}) + 3\sigma_{2}.(-12\sigma_{1}) + \\
    & \qquad  + 2\sigma_{1}.(18\sigma_{1}^{2} -12\sigma_{2})   + 3\sigma_{3}.(6) + 6\sigma_{1}^{3} -12\sigma_{1}.\sigma_{2} + 6\sigma_{3}
  \end{align*}
  and then
  \begin{align*}
 &  \Sigma_{3}N_{6} = \sigma_{1}^{3}.(120 - 144 + 36 - 72 + 18 + 36 + 6) + \\
 & \qquad +  \sigma_{1}.\sigma_{2}.(-144 -12 + 72 + 36 + 36 - 12 + 108 - 12 - 36 - 24 - 12) + \\
 & \qquad  +  \sigma_{3}.(36 + 12 - 36 - 36 + 18 + 6)  = 0
   \end{align*}

\parag{Conclusion} The system given by $\Sigma_{h}(F) = 0$ for each $h \in [2, k]$ which gives a characterization of the trace functions thanks to the previous results,  seems out of the range of an explicit computation  for $k$ big and we shall now try to find a computable family of elements in $W_{2}$ which gives a characterization of the trace functions.\\

The aim of the next section will be to construct such a nice explicit family which will be our candidate for such a characterization.

\section{Some second order PDE killing trace functions}

\subsection{Existence}

For $f$ an entire function on $\C$, the residue formula gives a integral formula for the trace of the  function $f$  computed from the symmetric functions $\sigma$ of $x_{1}, \dots, x_{k}$:
\begin{equation}
T(f)(\sigma) = \frac{1}{2i\pi}.\int_{\vert \zeta\vert = R} \ f(\zeta).\frac{P'_{\sigma}(\zeta).d\zeta}{P_{\sigma}(\zeta)} = \sum_{P_{\sigma}(x_{j}) = 0} f(x_{j}) 
\end{equation}
where $R$ is large enough compare to $\sigma$.\\
But for computing the second order partial derivatives of $T(f)$  this formula is not so convenient. \\
Remark that for $R$ large enough the quotient $P_{\sigma}(\zeta)\big/\zeta^{k}$ is near enough to $1$ on the circle $\vert \zeta \vert = R$, so we may integrate by part the previous formula and obtain a better formula for the computation of  partial derivatives in $\sigma$.

\begin{lemma}
For any entire function $f$ and any $\sigma \in \C^{k}$ we have, for $R$ large enough
\begin{equation}
T(f)(\sigma) = -\frac{1}{2i\pi}.\int_{\vert \zeta\vert = R} \ f'(\zeta).Log\big(P_{\sigma}(\zeta)\big/\zeta^{k}\big).d\zeta \  + \  k.f(0)
\end{equation}

\end{lemma} 

\parag{Proof} Thanks to the remark above, if we fix $\sigma$ in a compact subset of $N$, we may choose $R$ large enough to dispose of the holomorphic $1-$form $(Log P_{\sigma}(\zeta)\big/\zeta^{k}).d\zeta$ around the circle $\{ \vert \zeta \vert = R \}$ and to integrate by part the formula $(7)$. We have for $z$  near this circle :
\begin{equation*}
 \frac{d}{dz}\Big(f(z).Log\big(P_{\sigma}(z)\big/z^{k}\big)\Big) =  f'(z).Log\big(P_{\sigma}(z)\big/z^{k}\big) + f(z).P'_{\sigma}(z)\big/P_{\sigma}(z)  - k.f(z)\big/z
\end{equation*}
and this implies the formula $(8)$.$\hfill \blacksquare$\\

Now we may use this formula to compute the second order partial derivatives of $T(f)$. 

\begin{prop}\label{operators}
For any entire function $f$ on $\C$ the function $T(f)$ satisfies the following partial differential system, where we note $\partial_{h}:= \frac{\partial}{\partial \sigma_{h}}$:
\begin{align*}
& \big(\frac{\partial^{2}}{\partial_{p}\partial_{q}} - \frac{\partial^{2}}{\partial_{p+i}\partial_{q-i}}\big)[T(f)] = 0 \quad \forall p, q, i \quad {\rm such \ that} \quad p, q, p+i, q-i \in [1, k] \\
& \qquad \tag{@} \\
& \mathcal{T}^{m}[T(f)] = 0 \quad {\rm with} \quad \mathcal{T}^{m} := \partial_{1}.\partial_{m-1} + \big(\sum_{h=1}^{k} \sigma_{h}.\partial_{h}\big).\partial_{m} + \partial_{m}, \quad \forall m \in [2, k]
\end{align*}
\end{prop}

\parag{Proof}

\begin{align}
& (\partial_{p}[T(f)](\sigma) = (-1)^{p-1}.\frac{1}{2i\pi}.\int_{\vert \zeta \vert = R} f'(\zeta).\frac{\zeta^{k-p}.d\zeta}{P_{\sigma}(\zeta)} \\
& (\partial_{q}\partial_{p}[T(f)](\sigma) =  (-1)^{p+q}.\frac{1}{2i\pi}.\int_{\vert \zeta \vert = R} f'(\zeta).\frac{\zeta^{2k-q-p}.d\zeta}{P_{\sigma}(\zeta)^{2}}
\end{align}
and this depends only of the sum $p+q$, so this proves the first equations in $(@)$. \\
In the sequel we shall use the notation:
$$\mathcal{A}(p, q, i) :=  \big(\frac{\partial^{2}}{\partial_{p}\partial_{q}} - \frac{\partial^{2}}{\partial_{p+i}\partial_{q-i}}\big)$$
when $p, q, p+i, q-i$ are in $[1, k]$.\\

Let now $\mu$ be any integer between $[0, k-2]$ and write
\begin{equation*}
 \zeta^{\mu}.P_{\sigma}(\zeta) = \sum_{0}^{k} \ (-1)^{h}.\sigma_{h}.\zeta^{k-h+\mu}\quad {\rm with} \quad \sigma_{0}:= 1 .
\end{equation*}
Multiply by $f'(\zeta)$, divide by $P_{\sigma}(\zeta)^{2}$ and integrate on $\vert \zeta\vert = R$ leads to
\begin{equation*}
\partial_{k-\mu}[T(f)](\sigma)  =  \big(\sum_{h=0}^{k-1} \  \sigma_{h}.\partial_{k-\mu-1}\partial_{h+1} + \sigma_{k}.\partial_{k-\mu}.\partial_{k}\big)[T(f)](\sigma)
\end{equation*}
so the differential operator
\begin{equation}
\mathcal{T}_{0}^{\mu} :=  \sum_{h=0}^{k-1} \ \sigma_{h}.\partial_{k-\mu-1}\partial_{h+1} + \sigma_{k}.\partial_{k-\mu}.\partial_{k} + \partial_{k-\mu}
\end{equation}
kills the function $T(f)$ for any choice of the entire holomorphic function  $f: \C \to \C$.\\
Now, for $h \in [1, k-1]$ we have also that $\partial_{k-\mu-1}\partial_{h+1} - \partial_{k-\mu}.\partial_{h}$ also kills $T(f)$ (this differential operator  is equal to $-\mathcal{A}_{h, k-\mu, 1}$). So we obtain, with $m := k-\mu$
\begin{equation}
\mathcal{T}^{m} = \mathcal{T}_{0}^{k-m} - \sum_{h=1}^{k-1} \sigma_{h}.A_{h, m, 1} 
\end{equation}
also kills $T(f)$ for any $f$, proving the last equations in  $(@)$.$\hfill \blacksquare$\\

\parag{Remarks}\begin{enumerate}
\item All operators in the system $(@) $ have order $2$.
\item We have, for any $h \in [1, k]$:
\begin{equation}
[\partial_{h}, \mathcal{T}^{m}] = \partial_{m}.\partial_{h}
\end{equation}
\item In $(@)$ we may assume that $i$ takes only the value  $1$ because we have the relations
\begin{equation}
 \mathcal{A}_{p, q, i+1} =  \mathcal{A}_{p, q, i} + \mathcal{A}_{p+i, q-i, 1}
 \end{equation}
when $p, q, p+i+1, q-i-1$ are in $[1, k]$.
\item Using the formula 
$$ \frac{\partial x_{j}}{\partial \sigma_{h}} = (-1)^{h-1}.\frac{x_{j}^{k-h}}{P'_{\sigma}(x_{j})}\quad \forall h \in [1, k] $$
where $x_{j}$ is a local branch of root of $P_{\sigma}$ near a point $\sigma_{0} \in N \setminus \{\Delta(\sigma) = 0 \}$ one can check directly that the function $x_{j}^{q}$ is solution of the system $(@)$ for any $q \in \mathbb{N}$. But this computation is quite involved.\\
\end{enumerate}

\subsection{Symbols}

We shall study now the ideal in  $\C[\sigma_{1}, \dots, \sigma_{k,} \eta_{1}, \dots, \eta_{k}]$ (the algebra of polynomial functions on $T_{N}^{*}$)   generated by the symbols of the differential operators  in $(@)$  (note that $i = 1$ is enough from the remark 3. above).

\begin{lemma}\label{minors}
Let  $l_{\sigma}(\eta) := \sum_{h=1}^{k} \sigma_{h}.\eta_{h}$. Then the symbol of $\mathcal{T}^{m}$ is given by
\begin{equation}
Symb(\mathcal{T}^{m})(\sigma, \eta) = \eta_{1}.\eta_{m-1} + l_{\sigma}(\eta).\eta_{m}.
\end{equation}
and the ideal in $\C[\sigma, \eta]$ generated by the symbols of the operators $\mathcal{A}_{p, q, 1}$ and $\mathcal{T}^{m}$ is equal to the ideal generated by the $(2, 2)$ minors of the matrix
\begin{equation}
\begin{pmatrix} \eta_{1} & -l_{\sigma}(\eta) \\ \eta_{2} & \eta_{1} \\ \dots & \dots \\ \eta_{k} & \eta_{k-1}\end{pmatrix}
\end{equation}
\end{lemma}
\parag{Proof} It is obvious as the symbol of $\mathcal{T}^{m}$ is equal to $\eta_{1}.\eta_{m-1} + l_{\sigma}(\eta).\eta_{m}$.$\hfill \blacksquare$\\

One main ingredient in our proof of the theorem \ref{main th.} is the following proposition.

\begin{prop}\label{main prop.}
The ideal in $\C[\sigma_{1}, \dots, \sigma_{k,} \eta_{1}, \dots, \eta_{k}]$ generated by the $(2,2)-$minors of the matrix in $(16)$ is prime.
\end{prop}

\parag{Proof} The first step is to show that this ideal is reduced. This will be given by the proposition \ref{second step}. The second step will be given by the  lemma \ref{step 0}. $\hfill \blacksquare$\\

Denote by \ $\mathcal{I}_{Z}$ \ the ideal generated by the $(2,2)-$minors  $m_{\alpha}, \alpha \in A$ of the $(k, 2)-$matrix in $(16)$

\begin{lemma}\label{first step}
Each $\eta_{i}.\eta_{j}$ for $(i, j)$ in $ [1, k]\times [1, k]$ may be written as
\begin{equation}
 \eta_{i}.\eta_{j} = \sum_{\alpha }u_{i,j, \alpha}.m_{\alpha} + \eta_{k}.v_{i,j} 
 \end{equation}
where $u_{i,j,\alpha}$ are in $\C[\sigma_{1}, \dots, \sigma_{k}]$ (so  independent of $\eta$) and $v_{i,j}$ are polynomials in $(\sigma, \eta)$   linear in $\eta$ (so at most of degree $1$ in $\eta_{1}, \dots, \eta_{k}$).
\end{lemma}

\parag{Proof} We shall make a descending induction on $h \in [1, k]$ to prove the following assertion: \begin{itemize}
\item For a given $h \in [1, k]$ we can write $\eta_{j}.\eta_{h}$ as in $(17)$, for each $j \in [1, k]$.
\end{itemize}
As this is clearly true for $h = k$, let us assume  that for some $h \in [1, k-1]$  we have proved the assertion for $h+1$ and we shall prove it for $h$.\\
As $\eta_{1}.\eta_{h} + \sum_{p=1}^{k} \sigma_{p}.\eta_{p}.\eta_{h+1}$ is in $\mathcal{I}_{Z}$ the induction hypothesis gives  that $\eta_{1}.\eta_{h}$ may be written as in $(17)$.\\
For $j \in [2, k]$ we have $\eta_{j}.\eta_{h} - \eta_{j-1}.\eta_{h+1} \in \mathcal{I}_{Z}$ and then the induction hypothesis allows to completes the proof of our induction step. $\hfill \blacksquare$\\

\begin{lemma}\label{step 0}
Let $\vert Z\vert $ the algebraic subset  of common zeros of the polynomials in the ideal $\mathcal{I}_{Z}$.  Its Zariski open set  $\{ \eta_{k}\not= 0\}$ is dense, smooth and connected of dimension $k+1$. So $\vert Z\vert$ is irreducible.
\end{lemma}

\parag{Proof} Consider the polynomial map $\Phi : \C^{k}\times \C^{2} \to \C^{k}\times \C^{k}$ given by 
$$ \Phi(s, \zeta_{0}, \zeta_{1}) = (\sigma, \eta) \quad {\rm with} \quad  \sigma_{h} := (-1)^{h}.s_{h} \quad {\rm and} \quad \eta_{h}:= \zeta_{0}^{h}.\zeta_{1}^{k-h} \quad \forall  h \in [1, k].$$
Define also the hypersurface $\mathcal{H}$ in $\C^{k}\times \C^{2}$ 
$$ \mathcal{H} := \{(s, \zeta_{0}, \zeta_{1}) \in \C^{k}\times \C^{2} \ / \  \sum_{h=0}^{k} (-1)^{h}.s_{h}.\zeta_{0}^{h}.\zeta_{1}^{k-h} = 0 \}$$
with the convention $s_{0}\equiv 0$.\\
Remark that $H := \mathcal{H} \cap \{ \zeta_{0} = 1 \}$ is a smooth connected $k-$dimensional sub-manifold in $\mathcal{H}$ (isomorphic to $\C^{k}$ via the map $(s, \zeta_{1}) \mapsto (s_{1}, \dots, s_{k-1}, \zeta_{1})$).\\
Then we shall first verify that $\Phi(\mathcal{H}) \subset \vert Z \vert$. As the equation of $\mathcal{H}$ shows that on the image of $\mathcal{H}$ we have $l_{\sigma}(\eta) = -\zeta_{1}^{k}$ the vectors 
$$\begin{pmatrix} -l_{\sigma}(\eta) \\ \eta_{1}\\ \dots \\ \eta_{k-1}\end{pmatrix} = \zeta_{1}.\begin{pmatrix}\zeta_{1}^{k-1}\\ \zeta_{0}.\zeta_{1}^{k-2}\\ \dots \\ \zeta_{0}^{k-1}\end{pmatrix} \quad {\rm and} \quad \begin{pmatrix}\eta_{1}\\ \eta_{2}\\ \dots \\ \eta_{k}\end{pmatrix} = \zeta_{0}.\begin{pmatrix}\zeta_{1}^{k-1}\\ \zeta_{0}.\zeta_{1}^{k-2}\\ \dots \\ \zeta_{0}^{k-1}\end{pmatrix} $$
are proportional; this allows to conclude that $\Phi(\mathcal{H}) \subset \vert Z\vert$.\\
Now we have $\{\zeta_{0}\not= 0\} \cap \mathcal{H} \simeq H \times \C^{*}$ by the map $(s, \zeta_{0}, \zeta_{1}) \mapsto ((s, \zeta_{1}/\zeta_{0}), \zeta_{0})$ and $H$ is smooth of dimension $k$ and connected.\\
But on the open set $\{\zeta_{0}\not= 0\}\cap \mathcal{H}$ the restriction of the map $\Phi$ is an etale  $k-$sheeted covering onto $\{\eta_{k} \not= 0 \}\cap \vert Z\vert$:\\
In order to prove this we need  the following remark.
\parag{Remark} The open set $\{\eta_{k}\not= 0\} \cap \vert Z\vert$ is the complement of $\C^{k}\times \{0\}$ in $\vert Z\vert$. This means  that on $\vert Z\vert$ the equality $\eta_{k} = 0$ implies that $\eta = 0$ : \\
 looking of the $(k, 2)$ matrix whose $(2, 2)$ minors define the ideal  $\mathcal{I}_{Z}$ we immediately see that $\eta_{k} = 0$ implies $\eta_{k-1} = 0$ and then $\eta_{k-2}= \dots = \eta_{1} = 0$.\\
Of course the subset $\C^{k}\times \{0\}$ is in $\vert Z\vert$ because each equation of $\vert Z\vert$ is homogeneous of degree $2$ in $\eta$.\\

Now  on $\{\eta_{k} = 1\} \cap \vert Z\vert$ we see that the map $  \{\eta_{k} = 1\} \cap \vert Z\vert \to H$ given by $(\sigma, \eta) \mapsto (s, \eta_{k-1})$  is the inverse of   the restriction of $\Phi$ to $H \times \{\zeta_{0}= 1\}$; so it is an isomorphism. From this, it is easy to deduce, using homogeneity,  that the restriction of $\Phi$ to $\mathcal{H}\cap \{\zeta_{0}\not= 0\} \to \vert Z\vert \cap \{\eta_{k} \not= 0 \}$ is a $k-$sheeted etale cover corresponding to the $k-$th root of $\eta_{k}$.\\
 As $\mathcal{H}\cap \{\zeta_{0}\not= 0\}$ is isomorphic to $H \times \C^{*}$ which is smooth connected of dimension $k+1$, the only missing point is to prove the density of the open set $\{\eta_{k}\not= 0\} \cap \vert Z\vert$ in $\vert Z\vert$.\\
 So fix a point $(\sigma, 0) \in \C^{k}\times \{0\}$ and assume that $\sigma \not= 0$\footnote{This will be enough to prove our density statement !}. Let $z$ be a non zero root of the polynomial $\sum_{h=0}^{k} \sigma_{h}.z^{k-h} = 0$ and let $\varepsilon > 0$. Then, as the point $(s, \varepsilon, \varepsilon.z)$ is in $\mathcal{H}\cap \{\zeta_{0}\not= 0\}$, the point $\Phi(s, \varepsilon, \varepsilon.z)$ is in $\vert Z\vert \cap \{\eta_{k}\not= 0\}$ and when $\varepsilon \to 0$ it converges to $(\sigma, 0)$ concluding the proof.$\hfill \blacksquare$\\
 
 Note that the previous lemma implies that if  a holomorphic function $f$ on $N \times \C^{k}$ vanishes on $\vert Z\vert$, then $\Phi^{*}(f)$  vanishes on $\mathcal{H}$.

\begin{lemma}\label{stupid step}
Let $f$ be an holomorphic function on $N \times \C^{k}$ which is homogeneous of degree $0$ or $1$ in $\eta$ and which vanishes on $\vert Z\vert$. Then $f \equiv 0$.
\end{lemma}

\parag{Proof} In the case where $f$ is independent of $\eta$, then $\Phi^{*}(f)$ is independant of $\zeta_{0}$ and $\zeta_{1}$ and vanishes on $\mathcal{H}$. So $\Phi^{*}(f)$ vanishes. But $\Phi$ induces an isomorphism of $\C^{k}\times \{0\}$ to $\C^{k}\times \{0\}$ and then $f \equiv 0$.\\
In the case of a function linear in $\eta$  we may write $f = \sum_{h=1}^{k} f_{h}(\sigma).\eta_{h}$ and then $\Phi^{*}(f) = \sum_{h=1}^{k} g_{h}(s).\zeta_{0}^{h}.\zeta_{1}^{k-h}$. As $\Phi^{*}(f)$ vanishes on $\mathcal{H}$ we obtain\footnote{for instance put $\zeta_{0} = 1$ and use the fact that $1, z, \dots, z^{k-1}$ is a basis of $\mathcal{O}_{H}$ over $\mathcal{O}_{\C^{k}}$.} that $g_{h}\equiv 0$ for each $h \in [1, k]$ and then $f \equiv 0$.$\hfill \blacksquare$

\begin{prop}\label{second step}
Let $f$ be a holomorphic function on $N \times \C^{k}$ which vanishes on $\vert Z\vert$. Then $f$ is in $\mathcal{I}_{Z}$.
\end{prop}

\parag{Proof} It is enough to consider the case where $f$ is homogeneous in $\eta$ because the general case follows immediatly:\\
Consider the partial Taylor expansion  $f(\sigma, \eta) = \sum_{d=0}^{\infty} f_{d}(\sigma, \eta)$ where $f_{d}$ is homogeneous of degree $d$ in $\eta$. Then the homogeneity of $\vert Z\vert$ implies that $f$ vanishes on $\vert Z\vert$ if and only if each $f_{d}$ vanishes on $\vert Z\vert$.\\
Now, as the cases $d = 0, 1$ are already obtained in lemma \ref{stupid step}, consider the set of holomorphic functions on $N \times \C^{k}$ which are homogeneous of some degree $d \geq 2$ in $\eta$, vanish on $\vert Z\vert$ and are not in $\mathcal{I}_{Z}$. Assume that this set is not empty (if this is not the case, we are done, thanks to the lemma \ref{stupid step}). Assume that $f$ has minimal degree $d_{0} \geq 2$ in this set and write
$$ f(\sigma, \eta) = \eta_{k}.g(\sigma, \eta) + h(\sigma, \eta') $$
where $\eta' := (\eta_{1}, \dots, \eta_{k-1})$. Here $h$ is independent of $\eta_{k}$ homogeneous of degree $d_{0} \geq 2$ in $\eta'$ and $g$ is homogeneous of degree $d_{0}-1$ in $\eta$.\\
The function $h$ is in the ideal generated by $\eta_{i}.\eta_{j}, (i, j) \in [1, k-1]^{2}$ in $\mathcal{O}(\C^{k})[\eta']$ and, thanks to the lemma \ref{first step} we may write
$$ h(\sigma, \eta') = \sum_{\alpha} h_{\alpha}.m_{\alpha} + \eta_{k}.v(\sigma, \eta) $$
where $m_{\alpha}$ are the generators of $\mathcal{I}_{Z}$ (which are homogeneous of degree $2$ in $\eta$), where the functions $h_{\alpha}$ are holomorphic functions (independent of $\eta_{k}$) homogeneous of degree $d_{0}-2$ in $\eta'$ and where $v(\sigma, \eta)$ is homogeneous in $\eta$ of degree $d_{0}-1$. So we obtain
$$ f(\sigma, \eta) = \eta_{k}.(g - v)(\sigma, \eta) + \sum_{\alpha} h_{\alpha}.m_{\alpha} .$$
Now the function $\tilde{f} := f -  \sum_{\alpha} h_{\alpha}.m_{\alpha}$ is again vanishing on $\vert Z\vert$ (because the functions  $m_{\alpha}$ vanish on $\vert Z\vert$) and is homogeneous of degree $d_{0}$. As we have $\tilde{f} (\sigma, \eta) = \eta_{k}.\tilde{g}(\sigma, \eta)$ where $\tilde{g} := g - v$ is homogeneous of degree $d_{0}-1$ the irreducibility of $\vert Z\vert$ implies either $\eta_{k} = 0$ on $\vert Z\vert$ or $\tilde{g} = 0$ on $\vert Z\vert$. But we have seen in the lemma \ref{step 0}  that $\eta_{k}$ does not vanish at the generic point in $\vert Z\vert$. So we find a homogeneous degree $d_{0}-1$ function $\tilde{g}$ vanishing on $\vert Z\vert$ and the minimality of the degree $d_{0}$ of $f$ implies :
\begin{itemize}
\item either $d_{0}-1 \leq 1$ and then $\tilde{g} \equiv 0$ by the lemma \ref{stupid step} which implies $\tilde{f} \equiv 0$ and so  $f$ is in $\mathcal{I}_{Z}$. Contradiction.
\item or $d_{0}-1 \geq 2$ and then $\tilde{g}$ has to be in $\mathcal{I}_{Z}$ because of the minimality of $d_{0}$ and again $f$ is in $\mathcal{I}_{Z}$. Contradiction.
\end{itemize}
So every holomorphic function homogeneous in $\eta$ vanishing on $\vert Z\vert$ is in $\mathcal{I}_{Z}$.$\hfill \blacksquare$\\

\parag{Remark} It is an easy corollary of the previous proposition that the sheaf $\mathcal{I}_{Z}$ of ideals in $\mathcal{O}_{N \times \C^{k}}$ is reduced.$\hfill \square$\\

In the sequel we simply note $Z$ the common zero subset in $N \times \C^{k} \simeq T^{*}_{N}$ of the minors $(2, 2)$ of the matrix in $(16)$, which is the characteristic cycle (which is reduced) of the system $(@)$.\\

The simple corollary of the proposition \ref{main prop.} which follows will be used later on.

\begin{cor}\label{simplet 3}
The Zariski open set $Z \cap \{ \Delta(\sigma).\eta_{1}\not= 0 \}$ in   $Z$ is dense in $Z$.
\end{cor}

\parag{Proof} In fact we have seen that the set $\{\eta_{1} = 0 \} \cap Z$ is the union of $N\times \{0 \}$ and $\{\sigma_{k} = 0\} \cap Z$. But $Z$ has pure dimension $k+1$ and has a dimension $1$ fiber over each point in $N$. The  analytic subsets $N\times \{0 \}$, $\{\sigma_{k} = 0\} \cap Z$ and $\{\Delta(\sigma) = 0\} \cap Z$  have pure dimension $k$ so have no interior point in $Z$.$\hfill \blacksquare$\\

\begin{prop}\label{vanishing}
Let $\mathcal{U}$ be a non empty connected open set in $N$ and $Q$ a non zero  section on $\mathcal{U}$ of $D_{N}$ which annihilates any Newton polynomial $N_{m}(\sigma), m \in \mathbb{N}$ on $\mathcal{U}$. Then the symbol of $Q$ vanishes on $Z_{\vert \mathcal{U}} \subset (T_{N}^{*})_{\vert \mathcal{U}}$.
\end{prop}

\parag{Proof} It is enough to prove the proposition when $\mathcal{U}$ is the image by $s$ of the product  $U_{0}:= D_{1}\times \dots \dots D_{k}$ of $k$ disjoint discs in $\C$. In this case $\Xi$ induces an isomorphism on $\mathcal{U}$ and $Q = \Xi(P)$ where $P$ is a $\mathfrak{S}_{k}-$invariant differential operator on $U_{1}:= s^{-1}(\mathcal{U})$. Then we reduce the problem to show the following two facts :
\begin{enumerate}[i)]
\item If $P$ kills any Newton function $N_{m}(x), m \in \mathbb{N}$  on $U_{1}$ then the symbol of $P$ vanishes on $Y \cap (U_{1}\times \C^{k}) \simeq (T_{M}^{*})_{\vert U_{1}}$ where
$$ Y := \{(x, \xi_{1}, \dots, \xi_{k}) \in M\times \C^{k}  \ / \  \xi_{i}.\xi_{j}= 0 \quad \forall i \not= j \}.$$
\item The map $S\Xi^{*} : s^{-1}(T_{N}^{*}) \to T_{M}^{*}$ sends $s^{-1}(Z)$ into $Y$, and this induces a bijection between these sub-sets over the dense open set $s^{-1}\big(\{\Delta(\sigma) \not= 0 \}\big)$ in $M$.
\end{enumerate}
The proof of $i)$ will use the following two lemmas and the corollary \ref{vrac 3}.

 \begin{lemma}\label{step 2}
 Let $\mathcal{J}_{0}$ be the left ideal in $D_{M}$ generated by the differential operators  $\frac{\partial^{2}}{\partial x_{i}\partial x_{j}}$ for $(i, j) \in [1, k]^{2}, i \not= j$, and let $\mathcal{J}_{1}$ be the left ideal in $s_{*}(D_{M})^{\mathfrak{S}_{k}}$ which is given by its intersection with $s_{*}(\mathcal{J}_{0})$. Let $U_{1} \subset M$ be an open set as above and let $P$ be a $\mathfrak{S}_{k}-$invariant differential operator on $U_{1}$. 
Then $P$ may be written
 $$ P =  \sum_{q=0}^{N}\sum_{p=0}^{k-1} s^{*}(a_{p,q}).x_{j}^{p}.(\frac{\partial}{\partial x_{j}})^{q} \quad {\rm modulo} \quad \Gamma(\mathcal{U}, \mathcal{J}_{1}) $$
 where $a_{p, q}$ are holomorphic functions on $\mathcal{U} := s(U_{1})$.
 \end{lemma}
 
 \parag{Proof} Let $c.x^{\alpha}.(\frac{\partial}{\partial x})^{\beta}$ be a monomial of $P$ with $c \not= 0$. If $\beta \not= (0, \dots, 0, q, 0, \dots, 0)$ then the symmetrization $P_{1}$ of this monomial is in $\mathcal{J}_{1}$ and $P - P_{1} $ neither have this monomial nor these monomials deduced from it by the action of $\mathfrak{S}_{k}$. So we may write
 $$ P = \sum_{ \alpha}\sum_{q=0}^{N} c_{\alpha, j}.x^{\alpha}.(\frac{\partial}{\partial x_{j}})^{q} \quad {\rm modulo} \ \mathcal{J}_{1} .$$
 Now write $x^{\alpha} = (x(j))^{\alpha(j)}.x_{j}^{p}$, where we use the notation : $ x(j) := (x_{1}, \dots, \hat{x_{j}}, \dots, x_{k})$ and $\alpha(j) = (\alpha_{1}, \dots, \hat{\alpha_{j}}, \dots, \alpha_{k})$. The $\mathfrak{S}_{k}-$invariance of $P$ implies the invariance of $c_{\alpha, j}.(x(j))^{\alpha(j)}$  by the stabilizer $\mathfrak{S}_{k}(j)$ of $j$. Using the lemma
 \ref{simplet 2}, we obtain that $P$ has the following form
 $$ Q = \sum_{p= 0}^{k-1}\sum_{q = 0}^{N} \ s^{*}(a_{p, q}).x_{j}^{p}.(\frac{\partial}{\partial x_{j}})^{q}$$
 concluding the proof.$\hfill \blacksquare$\\
 
 \begin{lemma}\label{vrac 2}
 Let $x_{1}, \dots, x_{k}$ be distinct points in $\C$ and fix an integer $N$. Let $\mathcal{P}_{q}$ be the $\C-$vector space of polynomial of degree at most $q$ in $\C[z]$. Then the linear map
 $$ L : \mathcal{P}_{k.(N+1) -1} \to \C^{k.(N+1)} , \quad  Q \mapsto Q^{(p)}(x_{j}), \quad \forall  j \in [1, k], \quad \forall  p \in [0, N] $$
 is bijective.
 \end{lemma}
 
 \parag{Proof} This is, of course, a degenerate case of the standard Lagrange interpolation. We shall give a quick proof of  it  because the reader may be not so familiar with  this very degenerate case of Lagrange interpolation. Consider the ideal $\mathcal{I}$ in $\mathcal{O}_{\C}$ of germs vanishing at order $N$ at each point $x_{j}$ for $j \in [1, k]$. Then the exact sequence of coherent sheaves
 $$ 0 \to \mathcal{I} \to \mathcal{O} \to \mathcal{O}\big/\mathcal{I} \to 0 $$
 gives a surjective linear map $L_{1} :\mathcal{O}(\C) \to  \C^{k.(N+1)} $ analogous to $L$. Now each entire function $f$ on $\C$ may be written, thanks to Weierstrass division theorem\footnote{see [B-M 1] ch.II th. 3.2.9.}, in an unique way as 
 $$ f(z) = g(z).P(z)^{N+1} + q(z) $$
 where $g$ is an entire function, $P(z) := \prod_{j=1}^{k}(z - x_{j})$ and $q$ is a polynomial of degree at most $k.(N+1) - 1$. Then this allows to conclude that $L$ is surjective as $L_{1}(f) = L(q)$.$\hfill \blacksquare$\\
 
 \begin{cor}\label{vrac 3}
 Let $P$ be a $\mathfrak{S}_{k}-$invariant differential operator on $U_{1}$ which can be written as follows :
 $$P := \sum_{q=0}^{N}\sum_{p=0}^{k-1} s^{*}(a_{p,q}).x_{j}^{p}.(\frac{\partial}{\partial x_{j}})^{q} $$
where $a_{p, q}$ are holomorphic functions on \ $\mathcal{U} = s(U_{1})$. Assume that $P$ kills any Newton function $N_{m}(x), m \in \mathbb{N}$ on $U_{1}$. Then $P = 0$.
 \end{cor}
 
 \parag{Proof} It is enough to prove that each $a_{p,q}$ vanishes at each point in $\mathcal{U}\cap \{ \sigma_{k}\not= 0 \}$. Assume that  there exists $p_{0}, q_{0}$ and $\sigma^{0} = s(x_{1}^{0}, \dots, x^{0}_{k})  \in \mathcal{U}$ satisfying $\sigma_{k}^{0} \not= 0$ and such that  $a_{p_{0}, q_{0}}(\sigma^{0}) \not= 0$. Then we may find, thanks to the lemma \ref{vrac 2}, a polynomial $\Pi$ of degree at most $k.(N+1)-1$ such that $ \sum_{j=1}^{k} \Pi^{(q_{0})}(x^{0}_{j}).(x^{0}_{j})^{p_{0}} \not= 0$ and such that all derivatives at each point $x^{0}_{1}, \dots, x^{0}_{k}$ of order $\leq N$ and not equal to  $q_{0}$ (including $0$ when $q_{0}\not= 0$)  vanish. This is possible because for each $j \in [1, k]$ we have $x^{0}_{j}\not= 0$. Then we obtain a contradiction because $P$ has to kill the trace function
 $T(\Pi)(\sigma) := \sum_{j=1}^{k} \Pi(x_{j}) $ as $P$ has to kill each Newton function $N_{m}(x)$ and $T(\Pi)$ is a finite linear combination of these functions.
 But our choice of $\Pi$ implies $P[T(\Pi)] \not= 0$ at $\sigma^{0}$. Contradiction. So $P = 0$. $\hfill \blacksquare$\\
 
\parag{Proof of $i)$ in \ref{vanishing}} The three results above prove more that the fact $i) $ stated above because we obtain that $P$ is a section on $\mathcal{U}$  of the left ideal $\mathcal{J}_{1}$ of $s_{*}(D_{M})^{\mathfrak{S}_{k}}$. And the lemma \ref{simplet 0} shows that this ideal is globally generated by $S_{2}, \dots, S_{k}$ the elementary symmetric functions of the differential operators $\frac{\partial}{\partial x_{j}}, j \in [1, k]$.\\

The proof of $ii)$ in \ref{vanishing} will use the following two lemma.

\begin{lemma}\label{old-new}
Let $Z$ be the subspace of $\C^{k}\times \C^{k}$ defined by the ideal generated by the $(2, 2)$ minors of the matrix $(16)$. If $(\sigma, \eta)$ belongs to $Z$ and satisfies $\eta_{1}\not= 0$, then we have :
\begin{enumerate}
\item  $l_{\sigma}(\eta) \not= 0$;
\item  for each $h \in [1, k],\ \eta_{h} = \eta_{1}.(-\eta_{1}/l_{\sigma}(\eta))^{h-1}$;
\item  $l_{\sigma}(\eta)/\eta_{1}$ is a (non zero) root of the polynomial $P_{\sigma}$.
\end{enumerate}
\end{lemma}

\parag{Proof} The vanishing of the first minor equal to $\eta_{1}^{2} + l_{\sigma}(\eta).\eta_{2} $ implies that $l_{\sigma}(\eta)$ does not vanish. So 1. is proved.\\
Assume that we have proved that 2. is valid for $h \in [1, p]$ for $p \in [1, k-1]$. Note that this is clear for $h = 1$. Then the minor $\eta_{1}.\eta_{p} + \eta_{p+1}.l_{\sigma}(\eta)$ vanishes on $Z$ and gives
$$ \eta_{p+1} = -\eta_{1}.\eta_{p}/l_{\sigma}(\eta) = \eta_{1}.\big(-\eta_{1}/l_{\sigma}(\eta)\big)^{p} $$
using the induction hypothesis. So $2.$ is proved.\\
Now write $l_{\sigma}(\eta) = \sum_{h=1}^{k} \sigma_{h}.\eta_{h}$ and replace $\eta_{h}$ by the formula $2.$ This gives, after multiplication by $l_{\sigma}(\eta)^{k-1}$ and dividing by $\eta_{1}^{k}$
$$( l_{\sigma}(\eta)/\eta_{1})^{k} = \sum_{h=1}^{k} (-1)^{h-1}.\sigma_{h}.\big(l_{\sigma}(\eta)/\eta_{1}\big)^{k-h} $$
so $P_{\sigma}(\l_{\sigma}(\eta)/\eta_{1}) = 0$ proving $3.$.$\hfill \blacksquare$\\

In our next lemma we shall use now the notations introduced in the beginning of the paragraph 3.2.\\

\begin{lemma}\label{crucial}
Fix $\sigma \in \C^{k}$ and consider $\eta \in \C^{k}$ given by $\eta_{h} = a^{h-1}\eta_{1}, h \in [1, k]$ where $a \in \C^{*}$ is given. Then we have for $a.z \not= -1$:
$$ \sum_{h=1}^{k} \Theta_{h}(z, \sigma).\eta_{h} = \frac{(-a)^{k}}{1+a.z}.\big(P_{\sigma}(z) - P_{\sigma}(-1/a)\big).\eta_{1} $$
and  for $a.z = -1$:
$$   \sum_{h=1}^{k} \Theta_{h}(z, \sigma).\eta_{h} = z^{-k}.P'_{\sigma}(z).\eta_{1}$$
\end{lemma}

\parag{Proof} Note first that, by homogeneity of degree $1$ in $\eta$ we may assume that $\eta_{1} = 1$. As, by definition we have
$$ \sum_{h=1}^{k} \Theta_{h}(z, \sigma).\eta_{h} = \sum_{h=1}^{k} \Big(\sum_{p=0}^{h-1} (-z)^{h-p-1}.\sigma_{p}\Big).\eta_{h} $$
our hypothesis gives for $a.z \not= -1$:
\begin{align*}
&    \sum_{h=1}^{k} \Theta_{h}(z, \sigma).\eta_{h}  =  \sum_{h=1}^{k} (-a.z)^{h-1}\Big(\sum_{p=0}^{h-1} \sigma_{p}.(-z)^{-p}\Big) \\
&   \sum_{h=1}^{k} \Theta_{h}(z, \sigma).\eta_{h}   = \sum_{p=0}^{k-1} \sigma_{p}.(-z)^{-p}\Big(\sum_{h=p+1}^{k} (-a.z)^{h-1}\Big).\eta_{1} \\
&   \sum_{h=1}^{k} \Theta_{h}(z, \sigma).\eta_{h}   = \sum_{p=0}^{k-1} \sigma_{p}.(-z)^{-p}\Big((-a.z)^{p}\frac{(-a.z)^{k-p} -1}{-(1+a.z)}\Big)
\end{align*}
and then:
\begin{align*}
&   \sum_{h=1}^{k} \Theta_{h}(z, \sigma).\eta_{h}     =  -\frac{(-a)^{k}}{1+a.z}.\sum_{p=0}^{k-1} (-1)^{p}.\sigma_{p}.z^{k-p} + \frac{1}{1+ a.z} \sum_{p=0}^{k-1} \sigma_{p}.a^{p}\\
& {\rm but \  as} \quad (-a)^{k}.P_{\sigma}(-1/a) = (-a)^{k}.\sum_{p=0}^{k} (-1)^{p}.\sigma_{p}.(-1/a)^{k-p} = \sum_{p=0}^{k} \sigma_{p}.a^{p} \quad {\rm we \ obtain}\\
&   \sum_{h=1}^{k} \Theta_{h}(z, \sigma).\eta_{h}    =  -\frac{(-a)^{k}}{1+a.z}\Big( P_{\sigma}(z) - (-1)^{k}\sigma_{k}\Big) + \frac{(-a)^{k}}{1+a.z}\Big(P_{\sigma}(-1/a) - (-1)^{k}.\sigma_{k}\Big) \\ 
&   \sum_{h=1}^{k} \Theta_{h}(z, \sigma).\eta_{h}   =   -\frac{(-a)^{k}}{1+a.z}\Big( P_{\sigma}(z) - P_{\sigma}(-1/a)\Big)
\end{align*}

For $a.z = -1$ the computation gives, again for $\eta_{1} = 1$:
\begin{align*}
&   \sum_{h=1}^{k} \Theta_{h}(z, \sigma).\eta_{h}     =  \sum_{h=1}^{k} (-a.z)^{h-1}\Big(\sum_{p=0}^{h-1} \sigma_{p}.(-z)^{-p}\Big) \\
&   \sum_{h=1}^{k} \Theta_{h}(z, \sigma).\eta_{h}     =  \sum_{p=0}^{k-1} (k-p).\sigma_{p}(-z)^{-p} = z^{-k+1}.P'_{\sigma}(z)
\end{align*}
concluding the proof.$\hfill \blacksquare$\\

\parag{Remarks}\begin{enumerate}
\item When $z$ and $-1/a$ are distinct roots of $P_{\sigma}$ we find $ \sum_{h=1}^{k} \Theta_{h}(z, \sigma).\eta_{h} = 0$. If $z = -1/a$ is a double root of $P_{\sigma}$ we find again $ \sum_{h=1}^{k} \Theta_{h}(z, \sigma).\eta_{h}  = 0$. 
\item Assume now that  $\eta_{1} = \dots = \eta_{k-1} = 0$  and  $ z \not= 0$. Then we obtain
 $$\Theta(z, \sigma, \eta) = \frac{(-1)^{k-1}}{z}.\big(P_{\sigma}(z) - (-1)^{k}\sigma_{k}\big).\eta_{k}$$
 so  we find again $ \sum_{h=1}^{k} \Theta_{h}(z, \sigma).\eta_{h}  = 0$ when $\sigma_{k} = 0$ and  $P_{\sigma}(z) = 0$.
 \item At any point $(x, s(x), \eta)$  in $s^{-1}\big(Z \setminus (N \times \{0\})\big)$ with $\Delta(s(x)) \not= 0$ we have at least $(k-1)-$distinct  roots $x_{i}, i \in [1, k]$, of the polynomial $P_{s(x)}$ such that the numbers  $\sum_{h=1}^{k}\Theta_{h}(x_{i}, s(x)).\eta_{h}$ vanishes. So the symmetric functions of order $\geq 2$ of the numbers  $\sum_{h=1}^{k}\Theta_{h}(x_{i}, s(x)).\eta_{h}, i \in [1, k]$ vanish at such a point. 
\end{enumerate}

\parag{Proof of $ii)$ of \ref{vanishing}} As $Y$ is the union of the sub-manifolds
$$ Y_{j} := \{(x, \alpha_{1}, \dots, \alpha_{k}) \in T^{*}_{M} \ / \  \alpha_{1} = \dots = \hat{\alpha_{j}} = \dots = \alpha_{k} = 0 \}$$
for $j \in [1, k]$, it is enough to show that for any point $(x, s(x), \beta_{1}, \dots, \beta_{k})$  in $s^{-1}(Z)$ we have 
\begin{equation}
 \langle T_{s}^{*}(x, s(x), \beta_{1}, \dots, \beta_{k}), e_{j}\rangle = 0 
 \end{equation}
for $k-1$ values of $j \in [1, k]$. In fact it is enough to prove this fact assuming that $\beta_{1} \not= 0$ and $\Delta(s(x)) \not= 0$, because these two conditions define an open dense subset in $s^{-1}(Z)$ thanks to the corollary \ref{simplet 3}.

\bigskip

Then we obtain from lemma \ref{old-new} that there exists $j_{0}\in [1, k]$ such that we have
$$ \beta_{h} = \beta_{1}(-1/x_{j_{0}})^{h-1}, \ \forall h \in [1, k] \quad {\rm with} \quad P_{s(x)}(x_{j_{0}}) = 0.$$
Then the lemma \ref{crucial} with the remark 3. which follows it, gives $(18)$ for each $j \not= j_{0}$ concluding the proof of the proposition \ref{vanishing}.$\hfill \blacksquare$\\

\section{The solution}

\parag{Notations} Let $\mathcal{J}$ be the left ideal of $D_{N}$ generated by the elements in the system $(@)$.\\
Let $U_{0} := \sum_{h=1}^{k} h.\sigma_{h}.\partial_{h}$ and for $q \in \mathbb{N}$ let $\mathcal{K}_{q}$ the left ideal in $D_{N}$ generated by the differential operators in $(@)$ and $U_{0} -q$. So  $\mathcal{K}_{q} = \mathcal{J} + D_{N}.(U_{0} - q)$.

\subsection{The main theorem}

Our main result is the following characterization of ``trace functions''.

\begin{thm}\label{main th.}
Let $P \in \Gamma(N, D_{N})$ be a differential operator which  kills any Newton polynomial $N_{m}, m \in \mathbb{N}$. Then $P$ is in the left ideal $\Gamma(N, \mathcal{J})$ of $\Gamma(N, D_{N})$ 
\end{thm}

\parag{Proof} Remark first that $P$   kills any local  trace function on an open set $\mathcal{U}$ in $N$ by the lemma \ref{Newton} and so, thanks to the proposition \ref{vanishing}, if  $P$ is not $0$  the symbol of $P$ vanishes on $Z$. But using the proposition \ref{second step} we may find $P_{1} \in \Gamma(N,  \mathcal{J})$ such that the symbol of $P_{1}$ is equal to the symbol of $P$. Then $P - P_{1}$ again kills  any Newton polynomial $N_{m}, m \in \mathbb{N}$ and has order strictly less than the order of $P$. So if we assume that there exists such a $P  \in \Gamma(N, D_{N})$ which is not in $\Gamma(N, \mathcal{J})$, choosing such a $P$ with minimal order with these properties we obtain a contradiction by the previous argument. This conclude the proof.$\hfill \blacksquare$\\

\begin{cor}\label{surprise}
The $D_{N}-$module $\mathcal{M} :=  D_{N}\big/\mathcal{J}$ is sub-holonomic. For each $q \in \mathbb{N}$ the $D_{N}-$module $\mathcal{N}_{q} := D_{N}\big/\mathcal{K}_{q}$ is  holonomic on the complex manifold $N$. Its solutions of order $0$ on any simply connected open set \ $\mathcal{U}$ contained in
 $\{ \Delta(\sigma) \not= 0 \}$ are the linear combinations (with constant coefficients) of the $q-$th power of the branches on $\mathcal{U}$ of the multivalued function $z$ on $N$ which is the solution of the equation
 $$ P_{\sigma}(z) := z^{k} +  \sum_{h=1}^{k} (-1)^{h}.\sigma_{h}.z^{k-h} = 0 .$$
 \end{cor}
 
 \parag{Proof} As any $q-$th power of a local branch of the multivalued function $z$ is a local trace function, the differential operators in $(@)$ kills such a local branch. But any $q-$th power of such  a local branch is also killed by $U_{0} -q$.\\
 To see that $\mathcal{N}_{q}$ is holonomic it is enough to see that $Z \cap \{ \sum_{h=1}^{k} h.\sigma_{h}.\eta_{h} = 0 \}$ has dimension $k$. But this set is the union of $N \times \{0\}$ with the co-normal to the hyper-surface $\{\sigma_{k} = 0\}$ in $T_{N}^{*}$. $\hfill \blacksquare$\\
 
 \parag{Remark} It is easy to check the following commutation relations
 \begin{align*}
 &  A_{p, q}.U_{0} = \big(U_{0} - (p+q)\big).A_{p, q} \quad \forall p, q, p+1, q-1 \in [1, k] \\
 &  T^{m}.U_{0} = \big(U_{0} + m\big).T^{m} \quad \forall m \in [2, k]
 \end{align*}
 and this shows that  $\mathcal{J}.U_{0}$ is contained in $\mathcal{J}$. So we have an action of $U_{0}$ on the $D_{N}-$module $\mathcal{M}$ giving a spectral decomposition
 $$ \Gamma(N, \mathcal{M}) = \oplus_{m \in \mathbb{Z}} \ F_{m} $$
 where $F_{m}$ is the vector space of differential operators of pure weight $m$ modulo those in $\Gamma(N, \mathcal{J})$.\\
 Note that $P \in W_{2}$ has pure weight $p$ if and only if $[P, U_{0}] = -p.P$.  For instance the commutation relations above just mean that $A_{p, q}$ has pure weight $-(p+q)$ and that $T^{m}$ has pure weight $-m$.
 
 \begin{lemma}\label{derive}
 Let $\nabla := \sum_{h=0}^{k-1} (k-h).\sigma_{h}.\partial_{h+1}$. Then we have the following commutation relations
 \begin{align*}
 & [\nabla, \mathcal{T}^{h}] = -(k-h).\mathcal{T}^{h+1} \quad {\rm modulo} \quad (\partial_{1}.\partial_{h} -\partial_{2}.\partial_{h-1}) \quad \forall h \in  [2, k].\\
& [\nabla, A_{p, q, 1}] = - (k-p-1).A_{p+1, q, 1} - (k-q).A_{p, q+1} \\
&   {\rm for \, all} \  p, q \ {\rm such \ that} \  p, q, p+1, q-1, p+2, q-1 \ {\rm are \ in} \  [1, k]
\end{align*}
\end{lemma}

\parag{Proof} Recall that, by definition, we have $\mathcal{T}^{h} = \partial_{1}\partial_{h-1} + E.\partial_{h} + \partial_{h}$  for $h \in [2, k]$.
Note that we have $[\nabla, \partial_{h}] = - (k-h).\partial_{h+1}$ for each $h \in [1, k-1]$. So the second formula is  easy.\\
To prove the first formula, first remark that we have $[\sigma_{p}.\partial_{p+1}, E] = 0$ for each $p \in [1, k-1]$ and $[\partial_{1}, E] = \partial_{1}$,  so we have
\begin{align*}
& [\nabla, E] = k.\partial_{1} \quad {\rm and \ for} \  h \in [1, k-1] \\
& [\nabla, E.\partial_{h}] = \nabla. E.\partial_{h} - E.\partial_{h}.\nabla \\
&  [\nabla, E.\partial_{h}]  = (E.\nabla + k.\partial_{1}).\partial_{h} - E.(\nabla.\partial_{h} + [\partial_{h}, \nabla] ) \\
&  [\nabla, E.\partial_{h}] = k.\partial_{1}.\partial_{h} - (k-h).E.\partial_{h+1} \quad {\rm and} \  [\nabla, E\partial_{k}] = k.\partial_{1}\partial_{k}
\end{align*}

Then the first formula is now consequence of the following computation
\begin{align*}
&  [\nabla, \partial_{1}.\partial_{h-1}] = - (k-1).\partial_{2}.\partial_{h-1} - (k-h+1).\partial_{1}.\partial_{h}\\
& \qquad \qquad  = -(2k - h).\partial_{1}.\partial_{h} \quad {\rm modulo} \quad (\partial_{1}.\partial_{h} -\partial_{2}.\partial_{h-1}) 
\end{align*}
 for  each  $h \in [2, k-1] $ using our previous computations :
 \begin{align*}
 & [\nabla, \mathcal{T}^{h}] = [\nabla, \partial_{1}\partial_{h-1}] + [\nabla, E\partial_{h}] + [\nabla, \partial_{h}] \\
 & \qquad = -(2k-h).\partial_{1}\partial_{h} + k.\partial_{1}\partial_{h} - (k-h).\partial_{h+1} \  {\rm modulo} \ A_{1, h, 1} \\
 & \qquad = -(k-h).\mathcal{T}^{h+1}  \  {\rm modulo} \ A_{1, h, 1} 
 \end{align*}
 and we find $[\nabla, \mathcal{T}^{k} ] = 0 $ for $h = k$.$\hfill \blacksquare$\\
 
 So we have a right action of $\nabla$ on the $D_{N}-$module $\mathcal{M} = D_{N}\big/\mathcal{J}$.\\
 Note that $[\nabla, U_{0}] = \nabla$ as $\nabla$ is a derivation of weight $-1$. So the action of $\nabla$ on $\mathcal{M}$ shifts the weight decomposition of the global sections of $\mathcal{M}$ by $-1$.
 In the  bijective correspondence between global $0-$solutions\footnote{so $ \Gamma(N, Hom_{D_{N}}(\mathcal{M}, \mathcal{O}_{N}))$.} of $\mathcal{M}$ and $\Gamma(\C, \mathcal{O}_{\C})-$ given by $f \mapsto T(f)$ it  gives  the usual derivation of entire functions.\\
 
 Remark that, using the notation of the first sections, we have $\nabla = \Xi(S_{1}) = \Sigma_{1}$ because these derivations coincide on the elementary symmetric functions $\sigma_{1}, \dots, \sigma_{k}$.

\subsection{A holonomic system for trace forms}

For each $m \in \mathbb{Z}, m \geq -k+1$ and for each $\sigma \in N$ such that $\Delta(\sigma) \not= 0$ define
\begin{equation}
DN_{m}(\sigma) := \sum_{P_{\sigma}(x_{j}) = 0} \frac{x_{j}^{m+k-1}}{P'_{\sigma}(x_{j})} 
\end{equation}

\begin{prop}\label{DN.0}
Each $DN_{m}$ is the restriction to the open set $\{\Delta(\sigma) \not= 0 \}$ of a polynomial of (pure) weight $m$  in $\C[\sigma_{1}, \dots, \sigma_{k}]$ and the following properties are satisfied :
\begin{enumerate}[i)]
\item For $m \in [-k+1, -1], DN_{m} = 0 $.
\item For each $m \geq 1, \sum_{h=0}^{k} (-1)^{h}.\sigma_{h}.DN_{m-h} = 0 $ with the convention $\sigma_{0} \equiv 1$.
\item The polynomials $DN_{m}$ are in $\mathbb{Z}[\sigma_{1}, \dots, \sigma_{k}]$ and $DN_{1}, \dots, DN_{k}$ generate this $\mathbb{Z}-$algebra.
\item For $R \gg \vert\vert \sigma \vert\vert $ we have for each $m \geq -k+1$:
\begin{equation}
DN_{m}(\sigma) = \frac{1}{2i\pi}\int_{\vert \zeta\vert = R} \frac{\zeta^{m+k-1}.d\zeta}{P_{\sigma}(\zeta)} .
\end{equation}
\item For each $h \in [1, k]$ and  each $m \geq 0$ we have
\begin{equation}
\frac{\partial N_{m}}{\partial \sigma_{h}} = (-1)^{h-1}.m.DN_{m-h} .
\end{equation}
\end{enumerate}
\end{prop}

\parag{Proof} Remark first that the formula $(20)$ define a holomorphic function on $N$ because when $\sigma$ stays in a fixed  relatively compact set the integral is  independent of the choice of $R$ large enough. The residue formula implies easily that for $\Delta(\sigma) \not= 0$ this holomorphic function satisfies $(19)$.\\
Now, assuming that we have $\vert \sigma_{h}\vert \leq \varepsilon.R$ for each $h \in [1, k]$ with $\varepsilon \ll 1/k$ we obtain the estimate 
$$ \vert DN_{m}(\sigma) \vert \leq R^{m}\big(1+ k.\varepsilon\big).$$
 This implies that $DN_{m}$ is a polynomial of degree at most equal to $m$. The pure weight $m$ of this polynomial is easily obtained by a change of variable in $(20)$.This gives $i)$ and  the fact that $DN_{m}$ is in $\C[\sigma_{1}, \dots, \sigma_{k}]$ for each $ m \in \mathbb{N}$. \\
 With our definition of the polynomials $DN_{m}$ the formula $ii)$ is obvious. We show the first assertion in  $iii)$ by induction on $m \geq 0$ : \\
 For $m = 0$ the relation $\zeta^{k} = P_{\sigma}(\zeta) + \sum_{h=1}^{k} (-1)^{h-1}.\sigma_{h}.\zeta^{k-h}$ and the relation $ii)$ gives $DN_{0}\equiv 1$. Then the induction step is given by the relation $ii)$ as $\sigma_{0} = 1$. So $DN_{m}$ is in $\mathbb{Z}[\sigma_{1}, \dots, \sigma_{k}]$ for each $m \geq 0$.\\
 To complete the proof of $iii)$ it is enough to recall that for $h \in [1, k]$ the polynomial $DN_{h}$ is in $\mathbb{Z}[\sigma_{1}, \dots, \sigma_{h}]$, has pure weight $h$ and that the coefficient of $\sigma_{h}$ in it is equal to $(-1)^{h-1}$. \\
 The formula in  $v)$ is easily obtained by derivation  of the formula $(8)$ for $f(z) = z^{m}$ using $iv)$.$\hfill \blacksquare$\\
 
 \begin{defn}\label{DN.1} We shall call $DN_{m}$ the {\bf $m-$th derived Newton polynomial in $x_{1}, \dots, x_{k}$}.\\
 We shall say that $G \in \Gamma(N, \mathcal{O}_{N})$ is a {\bf trace form} if there exists $g \in \Gamma(\C, \mathcal{O}_{\C})$ such that
 \begin{equation}
 G(\sigma) = \sum_{j=1}^{k} \frac{g(x_{j})}{P'_{\sigma}(x_{j})}  := \tilde{T}(g)(\sigma) .
 \end{equation}
 \end{defn}
 
 \parag{Comment} Let $H := \{(\sigma, z) \in N \times \C \ / \ P_{\sigma}(z) = 0\}$ and $\pi : H \to N$ the projection. For $g \in \Gamma(\C, \mathcal{O}_{\C})$ we have the equality
 $$ Trace_{\pi}(g(z).dz\wedge d\sigma_{1}\wedge \dots\wedge d\sigma_{k-1}) = \tilde{T}(g)(\sigma).d\sigma_{1}\wedge \dots\wedge d\sigma_{k} $$
 which explains our terminology ``trace form''. This also prove the holomorphy of $G$ for any $g \in \Gamma(\C, \mathcal{O}_{\C})$.
 
 \begin{thm}\label{main for forms} Assume that the function $G \in \Gamma(N, \mathcal{O}_{N})$ is a trace form; then   it satisfies the system $(\tilde{@)}$ which is generated by the differential operators  $A_{p, q}$ for $p, q, p+1, q-1$ in $[1, k]$ and by  $\tilde{T}^{m} := T^{m} + \partial_{m}$ for $m \in [2, k]$. Moreover, for any $h \in [1, k] $ there exists a trace function $F_{h}$ such that $G = \partial_{h}F_{h}$.
 \end{thm}
 
 \parag{Proof} First remark that if $G$ is the trace form of the entire function $g = \sum_{m \geq 0} \gamma_{m}.z^{m}$, we have $G = \sum_{m \geq 0} \gamma_{m}.DN_{m-k+1}$ with uniform convergence on compact set. Thanks to formula $(21)$ the two assertions are easy consequences of the theorem \ref{main th.} using the following commutation relations
 \begin{align*}
 &  [A_{p, q}, \partial_{h}] = 0 \quad \forall p, q, p+1, q-1 \in [1, k] , \forall h \in [1, k] \\
 &  (T^{m} + \partial_{m}).\partial_{h} = \partial_{h}.T^{m} \quad \forall m \in [2, k],  \forall h \in [1, k]
 \end{align*}
 because they show that  when $F$ is a solution of $(@)$ then  $\partial_{h}F$ is solution of $(\tilde{@})$  for each $h \in [1, k]$.$\hfill \blacksquare$\\
 
 \parag{Remarks}\begin{enumerate}
 \item As an easy consequence of the previous theorem, using the fact that the characteristic variety of the left ideal $\tilde{\mathcal{J}}$ in $D_{N}$ generated by $(\tilde{@})$ is again equal to $Z$ which is irreducible, we obtain that the left ideal $\tilde{\mathcal{J}}$ in $D_{N}$ is the annihilator of all trace forms, or equivalently, of all polynomials $DN_{m}, m \in \mathbb{N}$.
 \item Any component of a global solution of the order $1$ linear system satisfied by the Lisbon integrals (see [B.MF 19]) is a trace form and then satisfies the system $(\tilde{@})$. Note that the components of such a solution are, up to the order and some signs,  the derivatives $\partial_{1}F, \dots, \partial_{k}F$ of a trace function $F$. So global solutions of the system described in {\it loc. cit.} may be identified with differential of  trace functions.
 \end{enumerate}

 \section{Miscellanous complements}

\subsection{Symmetric derivations}

\begin{lemma}\label{deriv.2}
The $\C[x_{1}, \dots, x_{k}]^{\mathfrak{S}_{k}}-$module of the symmetric derivations in  $W_{1}^{\mathfrak{S}_{k}}$ is generated by $U_{0}, \dots, U_{k-1}$ where we define for each $d \in \mathbb{N}$
$$ U_{d} := \sum_{j=1}^{k} \ \sum_{j=1}^{k} \ x_{j}^{d}.\frac{\partial}{\partial x_{j}} .$$
\end{lemma}

\parag{Proof} First remark that each $U_{d}$ for $d \geq k$ is in the $\C[x_{1}, \dots, x_{k}]^{\mathfrak{S}_{k}}-$module generated  by $U_{0}, \dots, U_{k-1}$ because, by definition we have for each integer $q$ the relation $ \sum_{h=0}^{k} (-1)^{h}.\sigma_{h}.x_{j}^{q+k-h} = 0$ for each $j \in [1, k]$, where $\sigma_{h} \in \C[x_{1}, \dots, x_{k}]^{\mathfrak{S}_{k}}$ is the $h-$th symmetric function of $x_{1}, \dots, x_{k}$.\\
Let now $\delta := \sum_{j=1}^{k} \ a_{j}(x).\frac{\partial}{\partial x_{j}}$ be a symmetric derivation. Writing for each $j \in [1, k]$  
$$a_{j}(x) := \sum_{p=0}^{N} a_{j,p}(x(j)).x_{j}^{p},$$
 where $x(j) := (x_{1}, \dots, \hat{x_{j}}, \dots, x_{k})$,  we see that the $\mathfrak{S}_{k}-$invariance of $\delta$ implies that for each $p$ the polynomial $a_{j,p}$ is $\mathfrak{S}_{k}[j]-$invariant, where $\mathfrak{S}_{k}[j]$ is the stabilizer of $j $ in $\mathfrak{S}_{k}$. Then, if $\sigma_{h}[j]$ denote the $h-$th symmetric function of $x(j)$, using the fact (see lemma \ref{simplet 2}) that $\sigma_{h}(j) = \sum_{q =0}^{h} \sigma_{h-q}.(-x_{j})^{q}$ for each $h \in [1, k-1]$ (with the convention $\sigma_{0}:= 1$), we may write
$$ \delta = \sum_{j=1}^{k} \ \sum_{p=0}^{N} \ b_{j,p}(\sigma_{1}, \dots, \sigma_{k}).x_{j}^{p}.\frac{\partial}{\partial x_{j}} = \sum_{p=0}^{N} \ b_{p}(\sigma_{1}, \dots, \sigma_{k}).U_{p} $$
because the $\mathfrak{S}_{k}-$invariance of $\delta$ implies that $b_{j,p}$ does not depend on $j \in [1, k]$ for each $p$.$\hfill \blacksquare$\\

\parag{Remark} As the derivation $U_{d}$ has pure weight $d-1$ this implies that the image by $\Xi$ of  any pure weight symmetric derivation in $W_{1}^{\mathfrak{S}_{k}}$ is an element of order $1$ in $W_{2}$ with pure weight at least equal to $-1$. As $\Xi$ induces an isomorphism on elements of order $0$, this implies that if  any derivation in $W_{2}$ which has pure weight $q$ is in the image of $\Gamma(N, \Xi) : W_{1}^{\mathfrak{S}_{k}} \to W_{2}$, then  $q \geq -1$. This shows that $\partial_{h}$ for $h \in [2, k]$ are not in this image. So $\Gamma(N, \Xi)$ is not surjective (and $\Xi$ also is not surjective).

\begin{lemma}\label{special} For each non negative integer $p$ define for $x_{1}, \dots, x_{k}$ with elementary symmetric functions $\sigma_{1}, \dots, \sigma_{k}$ such that $\Delta(\sigma) \not= 0$
\begin{equation}
\nabla_{p} := \sum_{j=1}^{p} \frac{x_{j}^{p}}{P'_{\sigma}(x_{j})}.\frac{\partial}{\partial x_{j}}
\end{equation}
Then, for any $Q \in \C[x_{1}, \dots, x_{k}]^{\mathfrak{S}_{k}}$, $\nabla_{p}(Q)$ is in $\C[x_{1}, \dots, x_{k}]^{\mathfrak{S}_{k}}$. So $\nabla_{p}$ define a derivation in $W_{2}$ with pure weight equal to $p-k$.\\
For each $p$ in $[0, k-1]$ we have
\begin{equation}
\nabla_{k-p} = (-1)^{p-1}.\frac{\partial}{\partial \sigma_{p}}
\end{equation}
\end{lemma}

\parag{Proof} The first assertion is an easy consequence of proposition \ref{DN.0} as we have
\begin{equation}
  \nabla_{h}[N_{m}] =  \sum_{j=1}^{k} \ (-1)^{h-1}.\frac{x_{j}^{k-h}}{P'_{\sigma}(x_{j})}.m.x_{j}^{m-1} = (-1)^{h-1}.m.DN_{m-h} \quad \forall p \in [1, k] \ \forall m \in \mathbb{N}
  \end{equation}

To prove $(24)$, as  a derivation in $W_{2}$ which vanishes on the Newton polynomials $N_{1}, \dots, N_{k}$ is nul, it is enough to compare the values of $\nabla_{k-p}$ and $\partial_{p}$ on $N_{m}$ for $m \in [1, k]$. But for any $m \in \mathbb{N}$ we have:
\begin{align*}
& \frac{\partial}{\partial \sigma_{h}}[N_{m}] = \sum_{j=1}^{k} \ m.x_{j}^{m-1}.\frac{\partial x_{j}}{\partial \sigma_{h}} = m.\sum_{j=1}^{k} (-1)^{h-1} \frac{x_{j}^{k-h+m-1}}{P'_{\sigma}(x_{j})} = (-1)^{h-1}.m.DN_{m-h} \\
\end{align*}
concluding the proof.$\hfill \blacksquare$\\

Note that the fact that the derivations $\nabla_{p}$ for $p \in [0, k-1]$ commute pair-wise is not obvious from a direct computation.\\
A consequence of the previous lemma is the fact that for any $P \in W_{2}$ of order $q$ then $\Delta^{2q}.P$ is in the image by $\Xi$ of $W_{1}^{\mathfrak{S}_{k}}$.\\

\subsection{Primitive Newton polynomials} 

For $m \geq k+1$ the differential form
\begin{equation}
\Omega_{m} := \sum_{h=1}^{k} (-1)^{h-1} \frac{N_{m-h}}{m-h}.d\sigma_{h} 
\end{equation}
is $d-$closed as a direct consequence of the formula $(21 )$. So it is $d-$exact. The following lemma shows that $ \Omega_{m} = d(PN_{m})$ for each $m \geq k+1$.

\begin{lemma}\label{PN}
Define for $m \geq k+1$
\begin{equation}
PN_{m}(\sigma) := \sum_{h=0}^{k} (-1)^{h-1}.\frac{N_{m-h}(\sigma)}{m-h}.\sigma_{h} 
\end{equation}
and for $m \in [1, k]$
\begin{equation}
PN_{m}(\sigma) := \sum_{h=0}^{m-1} (-1)^{h-1}.\frac{N_{m-h}(\sigma)}{m-h}.\sigma_{h} 
\end{equation}
with the convention $\sigma_{0} = 1$. Then we have for $m \geq 0$ :
\begin{equation}
\partial_{p}(PN_{m}) = (-1)^{m-p}.\frac{N_{m-p}}{m-p} \quad \forall p \in [1, k]
\end{equation}
with the following conventions : $N_{m-p} = 0$ for $m < p$ and $N_{m-p}/(m-p) = 1$ for $m = p$.
\end{lemma}

\parag{Proof} For $m \geq k+1$ the formula $(29 )$  is consequence of $(21)$ and the fact that $\sum_{h=0}^{k} (-1)^{h}.\sigma_{h}.DN_{m-p-h}(\sigma) = 0$. \\
For $m \leq k$ and $m \geq p+1$ the proof is the same, taking in account that $DN_{q} = 0$ for $q \in [-k+1, -1]$.\\
For $m < p$ the left hand-side is  clearly  $0$.  For $m = p$ the only term in $PN_{p}$ which depends on $\sigma_{p}$ is the term $(-1)^{p-1}.N_{p}\big/p$ and the coefficient of $\sigma_{p}$ in $N_{p}$ is equal to\footnote{Take the case $z^{p} - 1$ to compute it !} $(-1)^{p-1}.p$ concluding the proof.$\hfill \blacksquare$ \\

\parag{Example} we have for $k \geq 4$,  $PN_{1} =- \sigma_{1}, \quad PN_{2} =  (1/2)\sigma_{1}^{2} + \sigma_{2}$ \  and 
 $$ PN_{3} = \sigma_{3} - \sigma_{1}.\sigma_{2} + (1/6)\sigma_{1}^{3},$$
 $$ PN_{4} = -\sigma_{4} - \sigma_{1}.\sigma_{3} - (1/2)\sigma_{2}^{2} + (1/2)\sigma_{1}^{2}.\sigma_{2} - (1/12)\sigma_{1}^{4} \quad etc \dots$$
 and it is easy to check that $(29)$ is valid on these examples.

\parag{Final remarks} Using the commutation relations given in the proof of the theorem \ref{main for forms} we see that the left ideal in $\Gamma(N, D_{N})$ which annihilate all polynomial $PN_{m}$ for $m \in \mathbb{N}^{*}$ is generated by the differential operators $A_{p, q}, p, q, p+1, q-1 \in [1, k]$ and $T^{m}- \partial_{m}, m \in [2, k]$.\\
This system is given by differential operators with only order $2$ terms and has the same characteristic variety than $(@)$ or $(\tilde{@})$. So any $\sigma_{p}, p \in [1, k]$ is a solution of  this system.\\
Note that the polynomials $PN_{h}$ for $h \in [1, k]$ generate the algebra $\C[\sigma_{1}, \dots, \sigma_{k}]$ as this is true for the polynomials $N_{h}, h \in [1, k]$ .\\
It is easy to see that for $m \geq k+1$ the polynomial $PN_{m}(\sigma)$ is the trace  (via the map $\pi : H \to N$) of the restriction to $H$  of the polynomial $Q_{m, \sigma}[z]$ which is the primitive (in $z$) of the polynomial $z^{m-k-1}.P_{\sigma}(z)$ vanishing at $z = 0$.\\


\section{References}

\begin{itemize}
\item  {[B-M. 1]} Barlet, D. et Magnusson, J. {\it Cycles Analytiques Complexes I: Th\'eor\`{e}mes de pr\'eparation des cycles}, Cours Sp\'ecialis\'es $n^{0} 22$, Soci\'et\'e Math\'ematique de France (2014).
\item  {[B.MF. 19]} Barlet, D. and Monteiro Fernandes, T. {\it On Lisbon integrals}, matharxiv AG and CV  1906.09801
\item {[Bj]} \  Bjork, J.E.  {\it Rings of Differential Operators} North Holland (1979)
\item  {[Bor]} \ Borel, A. and al.  {\it Algebraic D-Modules}, Perspectives in Mathematics vol. 2 (1987) Academic Press.
\end{itemize}

\end{document}